\documentclass[reqno,12pt]{amsart}
\usepackage{amssymb, latexsym, euscript}
\newtheorem{theorem}{Theorem}[section]
\newtheorem{proposition}[theorem]{Proposition}

\newtheorem{lemma}[theorem]{Lemma}
\newtheorem{definition}[theorem]{Definition}
\newtheorem{corollary}[theorem]{Corollary}

\theoremstyle{definition}
\newtheorem{example}[theorem]{Example}
\newtheorem{remark}[theorem]{Remark}

\numberwithin{equation}{section}

   \def\sH{{\mathfrak H}}

      \def\dC{{\mathbb C}}

\def\cD{{\mathcal D}}      
   \def\cH{{\mathcal H}}   
   \def\cK{{\mathcal K}}   
      
      \def\cR{{\mathcal R}}

\def\bm\chi{\mbox{\boldmath$\chi$}}

\def\ker{{\rm ker\,}}

\def\dom{{\rm dom\,}}

\def\dim{{\rm dim\,}}
\def\diag{{\rm diag\,}}

\let\xker=\ker \def\ker{{\xker\,}}

\begin{document}
\title[On $J$-Self-Adjoint Operators with Stable $C$-Symmetry]{On $J$-Self-Adjoint Operators with Stable $C$-Symmetry}
\author[S.~Hassi]{Seppo~Hassi}
\author[S.~Kuzhel]{Sergii~Kuzhel}

\address{Department of Mathematics and Statistics \\
University of Vaasa \\
P.O. Box 700, 65101 Vaasa \\
Finland} \email{sha@uwasa.fi}

\address{Department of Applied Mathematics \\ AGH University of Science and Technology \\
30-059 Krakow, Poland}
 \email{kuzhel@mat.agh.edu.pl}



\begin{abstract}
The paper is devoted to a development of the theory of self-adjoint
operators in Krein spaces ($J$-self-adjoint operators) involving
some additional properties arising from the existence of
$C$-symmetries. The main attention is paid to the recent notion of
stable $C$-symmetry for $J$-self-adjoint extensions of a symmetric
operator $S$. The general results are specialized further by
studying in detail the case where $S$ has defect numbers $<2,2>$.
\end{abstract}

\maketitle

\section{Introduction}
Let $\mathfrak{H}$ be a Hilbert space with inner product
$(\cdot,\cdot)$ and let $J$ be a non-trivial fundamental symmetry,
i.e., $J=J^*$, $J^2=I$, and $J\not={\pm{I}}$.

The space $\mathfrak{H}$ equipped with the indefinite inner product
(indefinite metric) \ $[x,y]_J:=(J{x}, y), \
\forall{x,y}\in\mathfrak{H}$ is called  a Krein space
$(\mathfrak{H}, [\cdot,\cdot]_J)$.

An operator $A$ acting in $\mathfrak{H}$ is called
${J}$-self-adjoint if $A^*{J}={J}A$ i.e. if $A$ is
self-adjoint with respect to the indefinite metric
$[\cdot,\cdot]_J$.

The development of ${\mathcal{PT}}$-symmetric quantum mechanics (PTQM)
achieved during the last decade (see, e.g.,
\cite{B1,B2,MO,TT}) gives rise to a lot of new useful notions and interesting mathematical problems
in the theory of $J$-self-adjoint operators and, more generally, in the Krein space theory.

For instance, one of the key moments in PTQM
is the description of a hidden symmetry $C$ for a given pseudo-Hermitian Hamiltonian
in the sector of exact ${\mathcal{PT}}$-symmetry
\cite{B1,B2}. This immediately leads to the definition of $C$-symmetry for operators
acting in Krein spaces (Definition \ref{dad1}) and the problem of investigation of $J$-self-adjoint
operators with $C$-symmetry arises naturally.

In many cases, pseudo-Hermitian Hamiltonians of PTQM admit the
representation $A+V$, where a (fixed) self-adjoint operator $A$ and
a ${\mathcal{PT}}$-symmetric potential $V$ satisfy certain symmetry
properties which allow one to formalize the expression $A+V$ as a
family of $J$-self-adjoint\footnote{under a special choice of
fundamental symmetry $J$} operators $A_\varepsilon$ acting in a Krein space
$({\mathfrak H}, [\cdot,\cdot]_J)$. Here
$\varepsilon\in{\mathbb{R}}^m$ is a parameter characterizing the potential $V$.

One of important problems for the collection $\{A_\varepsilon\}$,
which is directly inspirited by PTQM, is the
description of quantitative and qualitative changes of spectra
$\sigma(A_\varepsilon)$.
Nowadays this topic has been analyzed with a wealth of technical tools (see, e.g.,
\cite{BKT,GRS,Z}).

If the potential $V$ is singular, then operators
$A_\varepsilon$ turn out to be $J$-self-adjoint extensions of the
\emph{symmetric} operator $S=A\upharpoonright\ker{V}$
and analysis of $A_\varepsilon$ can be carried out
by the extension theory methods \cite{A3,AK,AK1}.

In the present paper we continue such trend of investigation by studying the recent notion of stable $C$-symmetry for $J$-self-adjoint extensions of a symmetric operator $S$. This notion is natural
in the extension theory framework. Roughly speaking, if $A_\varepsilon$ belongs to the sector $\Sigma_J^{\textsf{st}}$ of stable $C$-symmetry (Definition \ref{dad2}), then $A_\varepsilon$ remains in $\Sigma_J^{\textsf{st}}$ under small variation of $\varepsilon$
(see, e.g., Theorem \ref{p14} with $\varepsilon=(\zeta, \phi, \xi, \omega)$).

The paper is structured as follows. Section 2 contains some preliminary results related to the Krein space theory and the
boundary triplets method in the extension theory. In our presentation, we have tried to emphasize the usefulness of the
Krein space ideology for the description of self-adjoint extensions of a symmetric operator.

In Section 3, $J$-self-adjoint extensions of a symmetric operator
$S$ with stable $C$-symmetry are investigated in a general setting.

Let $\{C_\alpha\}_{\alpha\in{D}}$ be the collection of operators $C$
(parameterized by a set of indexes ${D}$) which realize the property
of $C$-symmetry for $S$ (see Definition \ref{daad43}). The set
$\Upsilon$ of ${{J}}$-self-adjoint extensions $A\supset{S}$ which
commutes with any $C_\alpha$ plays a principal role in our
considerations. The description of $\Upsilon$ (Theorem \ref{neww33})
has some analogy with the description of self-adjoint extensions
$A\supset{S}$ which commute with a family of unitary operators
$\{U_\alpha\}$ satisfying the additional condition
$U\in\{U_\alpha\}\iff{U^*}\in\{U_\alpha\}$; see \cite{HK, KO, PH}.

If $\Upsilon$ is nonempty\footnote{the case $\Upsilon=\emptyset$ has
been considered in \cite{KVS}}, then there exists a boundary triplet
$(\mathcal{H}, \Gamma_0, \Gamma_1)$ of $S$ which provides the images
$\{\mathcal{C}_\alpha\}_{\alpha\in{D}}$ of
$\{C_\alpha\}_{\alpha\in{D}}$ in $\mathcal{H}$ with the preservation
of principal properties of $C_\alpha$ (Lemmas \ref{neww41},
\ref{neww67}). This enables one to describe $J$-self-adjoint
operators from the set $\Sigma_J^{\textsf{st}}$ of stable
$C$-symmetry in terms of stable $\mathcal{J}$-unitary operators of
the Krein space $(\mathcal{H}, [\cdot,\cdot]_{\mathcal{J}})$, where
$\mathcal{J}$ is the image of $J$ (Theorems \ref{neww80},
\ref{neww80b}). An `external' description of the family
$\{\mathcal{C}_\alpha\}_{\alpha\in{D}}$ (Theorem \ref{ttt6}) is
obtained by means of a reproducing kernel Hilbert space model
associated with the Weyl function of $S$ (cf. \cite{BHS}). This
leads to a characterization of the resolvents of $J$-self-adjoint
operators from the set $\Sigma_J^{\textsf{st}}$ of stable
$C$-symmetry (Theorem~\ref{absd1}), which is basic key to their
spectral analysis.

In Section 4, the set $\Sigma_J^{\textsf{st}}$ of stable
$C$-symmetry is studied under the assumptions that $S$ has defect
numbers $<2,2>$ and $S$ commutes with two anticommuting fundamental
symmetries $J$ and $R$ of the Hilbert space $\mathfrak{H}$ (see
(\ref{a2})). This condition means the existence of $J$-self-adjoint
extensions of $S$ with empty resolvent set \cite{KT} and it enables
one to investigate the spectral properties of
$A\in\Sigma_J^{\textsf{st}}$ in detail (Theorem \ref{wawa9}).
The key point here is the fact that the set
$\{C_\alpha\}_{\alpha\in{D}}$ of all $C$-symmetries of $S$ can be
expressed via $J$ and $R$ (Theorem \ref{uman4}).

As to the notations used in the paper: $\mathcal{D}(A)$ denotes the
domain of a linear operator $A$ and $A\upharpoonright{\mathcal{D}}$
means the restriction of $A$ onto a set $\mathcal{D}$. The symbols
$[A,B]:=AB-BA$ and $\{A,B\}:=AB+BA$ denote the commutator and
anti-commutator of the operators $A$ and $B$, respectively. The
symbols $\sigma(A)$ and $\rho(A)$ mean the spectrum and the
resolvent set of $A$.

\section{Preliminaries.}
\subsection{Elements of the Krein space theory.}
Let $(\mathfrak{H}, [\cdot,\cdot]_J)$ be a Krein space with
fundamental symmetry $J$. The corresponding orthogonal projections
$P_{\pm}=\frac{1}{2}(I{\pm}J)$ determine \emph{the fundamental
decomposition} of $\mathfrak{H}$:
\begin{equation}\label{d1}
\mathfrak{H}=\mathfrak{H}_+\oplus\mathfrak{H}_-, \qquad
\mathfrak{H}_-=P_{-}\mathfrak{H}, \quad \mathfrak{H}_+=P_{+}\mathfrak{H}.
\end{equation}

A subspace $\mathfrak{L}$ of $\mathfrak{H}$ is called \emph{hypermaximal neutral} if
$$
\mathfrak{L}=\mathfrak{L}^{[\bot]}=\{x\in\mathfrak{H}
 :  [x,y]_J=0, \ \forall{y}\in\mathfrak{L}\}.
$$

A subspace $\mathfrak{L}\subset\mathfrak{H}$ is called {\it uniformly positive
(uniformly negative)} if $[x,x]_J\geq{a}^2\|x\|^2$ \
($-[x,x]_J\geq{a}^2\|x\|^2$) $a\in\mathbb{R}$ for all
$x\in\mathfrak{L}$. The subspaces $\mathfrak{H}_{\pm}$ in (\ref{d1}) are
examples of uniformly positive and uniformly negative subspaces and
they possess the property of maximality in the corresponding classes
(i.e., $\mathfrak{H}_{+}$ ($\mathfrak{H}_{-}$) does not belong as a subspace to any
uniformly positive (negative) subspace).

Let $\mathfrak{L}_+(\not={\mathfrak H}_+)$ be a maximal uniformly
positive subspace. Then its $J$-orthogonal complement $
\mathfrak{L}_-=\mathfrak{L}_+^{[\bot]}$ is a maximal uniformly
negative subspace and the direct sum
\begin{equation}\label{d2}
\mathfrak{H}=\mathfrak{L}_+[\dot{+}]\mathfrak{L}_-
\end{equation}
gives another (in addition to (\ref{d1})) $J$-orthogonal
decomposition of $\mathfrak{H}$ into a
positive $\mathfrak{L}_+$ and a negative $\mathfrak{L}_-$ subspace of $\sH$. 


The subspaces $\mathfrak{L}_\pm$ in (\ref{d2}) can be described as
$\mathfrak{L}_+=(I+X)\mathfrak{H}_+$ and $\mathfrak{L}_-=(I+X^*)\mathfrak{H}_-$, where $X:\mathfrak{H}_+\to\mathfrak{H}_-$
is a strict contraction and $X^*:\mathfrak{H}_-\to\mathfrak{H}_+$ is the adjoint of $X$.

The self-adjoint operator $T=XP_++X^*P_-$ acting in $\mathfrak{H}$ is called
{\it a transition operator} from the fundamental decomposition
(\ref{d1}) to the decomposition (\ref{d2}). Obviously, $\mathfrak{L}_+=(I+T)\mathfrak{H}_+$
and $\mathfrak{L}_-=(I+T)\mathfrak{H}_-$.

Transition operators admit a simple description. Namely, a bounded operator $T$ in
$\mathfrak{H}$ is a transition operator if
and only if $T$ is a self-adjoint strict contraction and $\{J, T\}=0$.

The set of transition operators is in
one-to-one correspondence (via $\mathfrak{L}_\pm=(I+T)\mathfrak{H}_\pm$) with the set of decompositions
(\ref{d2}) of the Krein space $(\mathfrak{H},
[\cdot,\cdot]_J)$.

The projections $P_{\mathfrak{L}_\pm} : \mathfrak{H}\to\mathfrak{L}_{\pm}$ onto $\mathfrak{L}_\pm$ with
respect to the decomposition (\ref{d2}) are determined by the
formulas \cite{KK}:
$$
P_{\mathfrak{L}_-}=(I-T)^{-1}(P_--TP_+), \quad
P_{\mathfrak{L}_+}=(I-T)^{-1}(P_+-TP_-).
$$

The bounded operator
\begin{equation}\label{eae2}
{C}=P_{\mathfrak{L}_+}-P_{\mathfrak{L}_-}=J(I-T)(I+T)^{-1}
\end{equation}
also describes the subspaces ${\mathfrak{L}_\pm}$ in (\ref{d2}):
\begin{equation}\label{d4b}
\mathfrak{L}_+=\frac{1}{2}(I+{C})\mathfrak{H}, \qquad \mathfrak{L}_-=\frac{1}{2}(I-{C})\mathfrak{H}.
\end{equation}

The set of operators ${C}$ determined by (\ref{eae2}) is completely
characterized by the conditions
\begin{equation}\label{sos1}
{C}^2=I, \qquad  J{C}>0.
\end{equation}

It is known that the conditions in (\ref{sos1}) are equivalent to
the following representation of $C$ (see e.g. \cite[Lemma 2.8]{GKS}):
\begin{equation}\label{daad}
{C}=Je^{Y}, \qquad \{J, Y\}=0,
\end{equation}
where $Y$ is a bounded self-adjoint operator.

The decomposition $\mathfrak{H}=J\mathfrak{L}_+[\dot{+}]J\mathfrak{L}_-$ is called
\emph{dual} to (\ref{d2}). Its transition operator coincides with
$-T$ and the subspaces $J\mathfrak{L}_\pm$ of the dual decomposition are
described by (\ref{d4b}) with the adjoint operator ${C}^*$ instead
of ${C}$.

\subsection{Boundary triplets technique.}
Let $S$ be a symmetric (densely defined) operator with equal defect
numbers acting in a Hilbert space $\mathfrak{H}$.  A triplet
$(\mathcal{H}, \Gamma_0, \Gamma_1)$, where $\mathcal{H}$ is an
auxiliary Hilbert space and $\Gamma_0$, $\Gamma_1$ are linear
mappings of $\mathcal{D}(S^*)$ into $\mathcal{H}$, is called a
\emph{boundary triplet of} $S^*$ if the abstract Green identity
\begin{equation}\label{new2}
(S^*x, y)-(x, S^*y)=(\Gamma_1x, \Gamma_0y)_{\mathcal{H}}-(\Gamma_0x,
\Gamma_1y)_{\mathcal{H}}, \quad  x, y\in\mathcal{D}(S^*)
\end{equation}
is satisfied and the map $(\Gamma_0,
\Gamma_1):\mathcal{D}(S^*)\to\mathcal{H}\oplus\mathcal{H}$ is
surjective.

Let $(\mathcal{H}, \Gamma_0, \Gamma_1)$ be a boundary triplet for
$S^*$. Then the self-adjoint extensions
\begin{equation}\label{nnn1}
A_0=S^*\upharpoonright{\ker\Gamma_0}, \qquad \mbox{and} \qquad
A_1=S^*\upharpoonright{\ker\Gamma_1}
\end{equation}
are {\it transversal} in the sense that:
$$
\mathcal{D}(A_0)\cap\mathcal{D}(A_1)=\mathcal{D}(S) \qquad
\mbox{and}
 \qquad \mathcal{D}(A_0)+\mathcal{D}(A_1)=\mathcal{D}(S^*).
$$

These properties are characteristic for boundary triplets.
Precisely, for any transversal self-adjoint extensions $A_0$ and
$A_1$ of $S$ there exists a boundary triplet $(\mathcal{H},
\Gamma_0, \Gamma_1)$ of $S^*$ satisfying (\ref{nnn1}), see
\cite{DM}.

Fix a boundary triplet $(\mathcal{H}, \Gamma_0, \Gamma_1)$ for $S^*$
and consider the linear operators
\begin{equation}\label{new21}
\Omega_{+}=\frac{1}{\sqrt{2}}(\Gamma_1+i\Gamma_0), \qquad
\Omega_{-}=\frac{1}{\sqrt{2}}(\Gamma_1-i\Gamma_0)
\end{equation}
acting from $\mathcal{D}(S^*)$ into $\mathcal{H}$.
It follows from (\ref{new2}) and (\ref{new21}) that
\begin{equation}\label{new22}
(S^*x,y)-(x,S^*y)=i[(\Omega_{+}x,\Omega_{+}y)_{\mathcal
H}-(\Omega_{-}x,\Omega_{-}y)_{\mathcal H}].
\end{equation}
The formula (\ref{new22}) can be rewritten as:
\begin{equation}\label{new24}
(S^*x,y)-(x,S^*y)=i[\Psi{x},\Psi{y}]_{Z},
\end{equation}
where
\begin{equation}\label{new1}
\Psi=\left(\begin{array}{c} \Omega_+ \\
\Omega_- \end{array}\right) : \mathcal{D}(S^*) \to \textsf{H}=\left(\begin{array}{c} \mathcal{H} \\
\mathcal{H} \end{array}\right),
\end{equation}
maps $\mathcal{D}(S^*)$ into the Krein space $(\textsf{H},
[\cdot,\cdot]_{Z})$ with the indefinite metric
\begin{equation}\label{new30}
[\textsf{x},\textsf{y}]_{Z}=(x_0, y_0)-(x_1, y_1), \qquad \textsf{x}=\left(\begin{array}{c} x_0 \\
x_1\end{array}\right), \ \textsf{y}=\left(\begin{array}{c} y_0 \\
y_1\end{array}\right)\in\textsf{H},
\end{equation}
and the fundamental symmetry $Z=\diag(I,-I)$.

An arbitrary closed extension $A$ of $S$ is completely determined by
a subspace $\textsf{L}=\Psi\mathcal{D}(A)$ of $\textsf{H}$. In particular,
due to (\ref{new24}), $\Psi\mathcal{D}(A^*)=\textsf{L}^{[\perp]}$, where
$[\perp]$ means the orthogonal complement in the Krein space
$(\textsf{H}, [\cdot,\cdot]_Z)$. This leads to the following
statement:
\begin{lemma}\label{neww30}
Self-adjoint extensions of $S$ are in one-to-one correspondence with
hypermaximal neutral subspaces of the Krein space $(\textsf{H},
[\cdot,\cdot]_{Z})$.
\end{lemma}

The Weyl function $M(\cdot)$ and the characteristic function
$\Theta(\cdot)$ of $S$ associated with $(\mathcal{H}, \Gamma_0,
\Gamma_1)$ are defined as follows \cite{DM, Kochubei}:
\begin{equation}\label{neww65}
\begin{array}{c}
M(\mu)\Gamma_0f_\mu=\Gamma_1f_\mu, \quad
\forall{f}_\mu\in\ker(S^*-\mu{I}), \quad
\forall\mu\in\mathbb{C}_-\cup\mathbb{C_+}, \vspace{4mm} \\
\Theta(\mu)\Omega_{+}f_\mu=\Omega_{-}f_\mu,   \quad
\Theta(\overline{\mu})\Omega_{-}f_{\overline{\mu}}=\Omega_{+}f_{\overline{\mu}},
\quad \forall\mu\in\mathbb{C}_+.
\end{array}
\end{equation}
Notice, that $\Theta(\mu)$
and $M(\mu)$ are connected via the Cayley transform:
\begin{equation}\label{Cayley}
\Theta(\mu)=(M(\mu)-iI)(M(\mu)+iI)^{-1}, \quad \mu\in\mathbb{C}_+.
\end{equation}

\subsection{Fundamental decompositions depending on parameters.}
Let the decomposition (\ref{d2}) of a Krein space $(\mathfrak{H},
[\cdot,\cdot]_J)$ depend on parameters:
\begin{equation}\label{d2c}
\mathfrak{H}=\mathfrak{L}_+^\alpha[\dot{+}]\mathfrak{L}_-^\alpha, \qquad
\alpha=(\alpha_1,\ldots,\alpha_m)\in{D}\subset{\mathbb{R}^m},
\end{equation}
where $\alpha$ runs a domain ${D}\subset{\mathbb{R}^m}$. We are
going to illustrate what may happen with decomposition (\ref{d2c})
when $\alpha$ tends to a certain point $\alpha_0$ lying on the
boundary $\partial{D}$ of $D$. Observe, that the subspaces
$\mathfrak{L}_\pm^\alpha$ are determined uniquely by the family of
transition operators $\{T_\alpha\}_{\alpha\in{D}}$.

Assume that $T_\alpha$ tends (in the strong sense) to a unitary
operator $Q$ in $\mathfrak{H}$ when $\alpha\to\alpha_0$. Then $Q$ is
a unitary and self-adjoint operator such that $\{J,Q\}=0$. In that
case, elements $x_++T_\alpha{x_+}$ \ $({x}_+\in\mathfrak{H}_+)$ of
$\mathfrak{L}_+^\alpha$ converge to elements of the subspace \ $
\mathfrak{L}=\{ x_+ + Qx_+ :\, {x_+}\in\mathfrak{H}_+ \}. $

On the other hand, elements $x_-+T_\alpha{x_-}$ \ $({x}_-\in\mathfrak{H}_-)$ of
$\mathfrak{L}_-^\alpha$ converge to elements $x_-+Q{x_-}$ of the \emph{same} subspace $\mathfrak{L}$
(since $x_-+Q{x_-}=x_++Q{x_+}$ for $x_-=Qx_+$). This means that
the `pointwise limit' of $\mathfrak{L}_\pm^\alpha$ (as $\alpha\to\alpha_0$) \emph{coincides with the
hypermaximal neutral subspace}
$$
\mathfrak{L}=(I+Q)\mathfrak{H}_+=(I+Q)\mathfrak{H}_-=(I+Q)\mathfrak{H}
$$
of the Krein space $(\mathfrak{H},
[\cdot,\cdot]_J)$.

Similarly, the subspaces
$J\mathfrak{L}_\pm^\alpha$ of the dual decomposition
\begin{equation}\label{d2d}
\mathfrak{H}=J\mathfrak{L}_+^\alpha[\dot{+}]{J}\mathfrak{L}_-^\alpha
\end{equation}
`tend' to the dual hypermaximal neutral subspace (as $\alpha\to\alpha_0$)
$$
\mathfrak{L}^\sharp=J\mathfrak{L}=(I-Q)\mathfrak{H}_+=(I-Q)\mathfrak{H}_-=(I-Q)\mathfrak{H}.
$$
Therefore, the limits ($\alpha\to\alpha_0$) of the decompositions
(\ref{d2c}) and (\ref{d2d}) give rise to a new decomposition
\begin{equation}\label{sas3}
\mathfrak{H}=\mathfrak{L}[\dot{+}]\mathfrak{L}^\sharp,
\end{equation}
which is fundamental in the \emph{new Krein space} $(\mathfrak{H},
[\cdot,\cdot]_Q)$ with the fundamental symmetry $Q$.

The following examples illustrate the phenomena described above.

\begin{example}\label{ex1}
Let $R$ be a unitary and self-adjoint operator
in $\mathfrak{H}$ which anti-commutes with $J$, i.e., $\{J,R\}=0$.
Then the operator
\begin{equation}\label{neww12}
R_{\omega}={R}e^{i{\omega}J}={R}[\cos\omega+i(\sin\omega)J], \quad
\omega\in[0,2\pi),
\end{equation}
is also unitary and self-adjoint in $\mathfrak{H}$ and, furthermore, $\{J,
R_\omega\}=0$. According to (\ref{daad}), the family of operators
\begin{equation}\label{sas44}
{C}_{\chi,\omega}=Je^{\chi{R_\omega}}=J[(\cosh\chi)I+(\sinh\chi)R_\omega],
\quad \chi\in\mathbb{R},
\end{equation}
is defined by (\ref{eae2}), and the corresponding family of
transition operators $T_{\chi,\omega}$ has the form
$$
T_{\chi,\omega}=-(\tanh\frac{\chi}{2})R_\omega,  \quad
\chi\in\mathbb{R}.
$$
The operators ${{C}}_{\chi, \omega}$ (or $T_{\chi,\omega}$) determine the
subspaces
$$
\mathfrak{L}_\pm^{\chi,\omega}=\frac{1}{2}(I\pm{{C}_{\chi,\omega}})\mathfrak{H}=(I+T_{\chi,\omega})\mathfrak{H}_{\pm},
\quad \alpha=(\chi, \omega)\in{D}=\mathbb{R}\times[0,2\pi).
$$
If $\chi\to\infty$, then $T_{\chi,\omega}$ tends to a unitary
operator $Q=-R_\omega$. The corresponding subspaces
$\mathfrak{L}_\pm^{\chi,\omega}$ converge pointwise to the
hypermaximal neutral subspace
$\mathfrak{L}=(I-R_\omega)\mathfrak{H}$ of the Krein space
$(\mathfrak{H},[\cdot,\cdot]_J)$. Similarly, when $\chi\to\infty$,
the subspaces of the dual decomposition
$$
J\mathfrak{L}_\pm^{\chi,\omega}=\frac{J}{2}(I\pm{{C}_{\chi,\omega}})\mathfrak{H}=
\frac{1}{2}(I\pm{{C}_{-\chi,\omega}})J\mathfrak{H}=\frac{1}{2}(I\pm{{C}_{-\chi,\omega}})\mathfrak{H}
$$
tend to the dual hypermaximal neutral subspace
$\mathfrak{L}^\sharp=(I+R_\omega)\mathfrak{H}$. The 'limiting'
subspaces $\mathfrak{L}$ and $\mathfrak{L}^\sharp$ give rise to the
fundamental decomposition of the Krein space
$(\mathfrak{H},[\cdot,\cdot]_{Q})$ with $Q=-R_\omega$.
\end{example}


\begin{example}\label{ex2}
Let $(\mathcal{H}, \Gamma_0, \Gamma_1)$ be a
boundary triplet of a symmetric operator $S$ and let $(\textsf{H},
[\cdot,\cdot]_{Z})$ be the corresponding Krein space defined by
(\ref{neww30}). The fundamental decomposition of $(\textsf{H},
[\cdot,\cdot]_{Z})$ has the form
\begin{equation}\label{new25}
\textsf{H}=\left(\begin{array}{c} \mathcal{H} \\
0 \end{array}\right)\oplus\left(\begin{array}{c} 0 \\
\mathcal{H}
\end{array}\right), \qquad \textsf{H}_+=\left(\begin{array}{c} \mathcal{H} \\
0 \end{array}\right), \quad \textsf{H}_-=\left(\begin{array}{c} 0 \\
\mathcal{H}
\end{array}\right).
\end{equation}
According to (\ref{new24}), the subspaces
\begin{equation}\label{new91}
\textsf{L}_+^\mu=\Psi\ker(S^*-\mu{I}), \quad
\textsf{L}_-^\mu=\Psi\ker(S^*-\overline{\mu}{I}), \quad
\mu\in\mathbb{C}_+
\end{equation}
are, respectively, uniformly positive and uniformly negative in the
Krein space $(\textsf{H}, [\cdot,\cdot]_{Z})$ and they form a
$Z$-orthogonal decomposition
\begin{equation}\label{d2cd}
\textsf{H}=\textsf{L}_+^\mu[+]\textsf{L}_-^\mu \qquad
\mu\in{D=\mathbb{C}_+},
\end{equation}
cf. (\ref{d2c}). The corresponding family of transition operators
$\{T_\mu\}_{\mu\in{\mathbb{C}_+}}$ from the fundamental
decomposition (\ref{new25}) to (\ref{d2cd}) has the operator-matrix
form (with respect to (\ref{new25})):
\begin{equation}\label{nnn2}
T_\mu=\left(\begin{array}{cc} 0 & \Theta(\overline{\mu}) \\
\Theta(\mu) & 0
\end{array}\right) \qquad \mu\in{\mathbb{C}_+}, \qquad
\Theta(\overline{\mu})=\Theta^*(\mu),
\end{equation}
where $\Theta(\mu)$ is the characteristic function of $S$ associated
with the boundary triplet $(\mathcal{H}, \Gamma_0, \Gamma_1)$.
\end{example}

\begin{remark}\label{vvv1}
The decomposition (\ref{d2cd}) depends on the choice of the boundary
triplet $(\mathcal{H}, \Gamma_0, \Gamma_1)$ which is determined by
the mapping $\Psi$ (see (\ref{new1})):
\begin{equation}\label{www2}
\left(\begin{array}{c} {\Gamma}_0x \\
{\Gamma}_1x \end{array}\right)=B{\Psi}x, \quad
x\in\mathcal{D}(S^*), \quad B=\frac{1}{\sqrt{2}}\left(\begin{array}{cc} -iI & iI \\
I & I \end{array}\right).
\end{equation}

Let $U$ be an arbitrary $Z$-unitary operator in the Krein space
$(\textsf{H}, [\cdot,\cdot]_{Z})$ and let ${\Psi}_U=U\Psi$.
Taking (\ref{www2}) into account,
we conclude that $U$ determines a new boundary triplet $(\mathcal{H},
{{\Gamma}}_0^U, {\Gamma}_1^U)$ of $S^*$, where
$$
\left(\begin{array}{c} {\Gamma}_0^Ux \\
{\Gamma}_1^Ux \end{array}\right)=B{\Psi}_Ux, \quad
x\in\mathcal{D}(S^*).
$$
In that case, $\textsf{H}=U\textsf{L}_+^\mu[+]U\textsf{L}_-^\mu$
($\textsf{L}_\pm^\mu$ are defined by (\ref{new91})); the family of
transition operators $\{T_\mu^U\}_{\mu\in{\mathbb{C}_+}}$ is
determined in the same manner as (\ref{nnn2}) by the characteristic
function $\Theta_U(\mu)$ of $S$ in $(\mathcal{H}, {{\Gamma}}_0^U,
{\Gamma}_1^U)$ and $\Theta_U(\mu)$ can be expressed via the
Krein-Shmul'yan transformation (see \cite{DHMS1}):
$$
\Theta_U(\mu)=(U_{10}+U_{11}\Theta(\mu))(U_{00}+U_{01}\Theta(\mu))^{-1},
$$
where the bounded operators $U_{ij}\in\mathcal{B}[\mathcal{H}]$
originate from the decomposition $U=(U_{ij})_{i,j=0}^1$ with respect
to (\ref{new25}).
\end{remark}

Let $r$ be a real point of regular type of $S$. Then the
operator
\begin{equation}\label{sas99}
A_{r}=S^*\upharpoonright{{\mathcal{D}}(A_{r})}, \qquad
{\mathcal{D}}(A_{r})={\mathcal{D}}(S)\dot{+}\ker(S^*-rI)
\end{equation}
is a self-adjoint extension of $S$ and
$\textsf{L}_{r}=\Psi\mathcal{D}(A_{r})$ is a hypermaximal neutral subspace in
$(\textsf{H}, [{\cdot}, \cdot]_{Z})$. Therefore, $\textsf{L}_{r}=(I+X_{r})\textsf{H}_+$,
where $X_{r}$ is a \emph{unitary} mapping of $\textsf{H}_+$ onto
$\textsf{H}_-$.

It follows from the results of \cite{Kochubei, SH} that
$$
X_{r}=\Theta(r)=s-\lim_{\mu\to{{r}}}\Theta(\mu),
 \quad X_{r}^{-1}=X_{r}^{*}=\Theta^{-1}(r)=s-\lim_{\overline{\mu}\to{{{r}}}}\Theta(\overline{\mu}),
$$
where $\mu\in\mathbb{C}_+$.
This means that the transition operators
$\{T_\mu\}_{\mu\in{\mathbb{C}_+}}$ determined by (\ref{nnn2}) converge
(in the strong sense) to
$$
Q_{r}=\left(\begin{array}{cc}
0 & \Theta^{-1}(r) \\
\Theta(r) & 0 \end{array}\right).
$$
The operator $Q_{r}$ is self-adjoint and unitary in the Hilbert space $\textsf{H}$ (with the
inner product $(\cdot,\cdot)=[Z\cdot,\cdot]_Z$) and $\{Z,Q_{r}\}=0$.

Thus, if  $\mu\in\mathbb{C}_+$ tends to a real point
 of regular type $r$ of $S$, then the subspaces $\textsf{{L}}_+^\mu$ and
$\textsf{{L}}_-^{{\mu}}$ in (\ref{d2cd}) converge to the hypermaximal neutral subspace
$\textsf{L}_{r}=(I+\Theta(r))\textsf{H}_+=(I+\Theta^{-1}(r))\textsf{H}_-=(I+Q_{r})\textsf{H}$.

Furthermore, the dual subspaces $Z\textsf{{L}}_\pm^\mu$ tend to the dual
hypermaximal neutral subspace
$\textsf{L}_{r}^\sharp=Z\textsf{L}_{r}=(I-\Theta(r))\textsf{H}_+=(I-\Theta^{-1}(r))\textsf{H}_-=(I-Q_{r})\textsf{H}$
of the Krein space $(\textsf{H}, [\cdot, \cdot]_Z)$. The subspaces
$\textsf{L}_{r}$ and $\textsf{L}_{r}^\sharp$ give rise to a decomposition
$\textsf{H}=\textsf{L}_{r}[\dot{+}]\textsf{L}^\sharp_{r}$, which is
fundamental in the new Krein space $(\textsf{H},
[\cdot,\cdot]_{Q_{r}})$.

%

\section{$J$-self-adjoint operators with stable ${C}$-symmetry.}
\subsection{Definition of $C$-symmetry and stable ${C}$-symmetry.}
By analogy with \cite{B1, B2} the definition of
${C}$-symmetry in the Krein space setting can be formalized as follows.
\begin{definition}\label{dad1}
An operator ${A}$ in a Krein space $(\mathfrak{H},
[\cdot,\cdot]_{{J}})$ has the property of ${C}$-symmetry if
there exists a bounded linear operator ${C}$ in $\mathfrak{H}$ such
that: \ $(i) \ {C}^2=I;$ \quad $(ii) \ {J}C>0$; \quad
$(iii) \ A{C}={C}A$.
\end{definition}

In what follows, in order to avoid possible misunderstandings, we
will use the sentence:  \emph{an operator $A$ has the  property of
$C$-symmetry} in the meaning of the general property (realized by
some (non-specified) operator $C$). On the other hand, the sentence:
\emph{an operator $A$ has the $C$-symmetry} will mean that the
property of $C$-symmetry of $A$ is realized by the concrete operator
$C$.

In view of (\ref{eae2}) and (\ref{sos1}) the property of
${C}$-symmetry of $A$ means that $A$ can be decomposed with respect
to the decomposition $(\ref{d2})$ (with subspaces $\mathfrak{L}_\pm$
determined by (\ref{d4b})):
\begin{equation}\label{neww4}
A=A_+\dot{+}A_-, \qquad
A_+=A\upharpoonright{\mathfrak{L}_+}, \quad
A_-=A\upharpoonright{\mathfrak{L}_-}.
\end{equation}

If a ${J}$-self-adjoint\footnote{see the definition in Introduction.} operator $A$ possesses the property
of $C$-symmetry, then its components $A_\pm$ in
(\ref{neww4}) turn out to be self-adjoint operators in the Hilbert
spaces ${\mathfrak{L}_+}$ and ${\mathfrak{L}_-}$ with the inner
products ${[\cdot,\cdot]_{J}}$ and
$-{[\cdot,\cdot]_{J}}$, respectively. This simple
observation leads to the following statement.

\begin{proposition}[\cite{K(Krein)}]\label{sese1}
A ${J}$-self-adjoint operator $A$ has the
property of ${C}$-symmetry if and only if $A$ is similar to a
self-adjoint operator in $\mathfrak{H}$. If a
${J}$-self-adjoint operator $A$ has the ${C}$-symmetry then
the adjoint operator $A^*$ has the ${C}^*$-symmetry.
\end{proposition}

\begin{remark} The first part of Proposition \ref{sese1}
(\emph{property of $C$-symmetry $\Leftrightarrow$ similarity to a self-adjoint one}) is a
direct consequence of the Phillips theorem \cite[Theorem 6.1]{PH},
which relates `the similarity' of $A$ and `the invariance' of $A$
with respect to (\ref{d2}) to each other. In this context, the
concept of $C$-symmetry allows one to specify the decomposition (\ref{d2})
and to describe a self-adjoint operator mentioned in Proposition \ref{sese1} explicitly. Indeed, if a
${J}$-self-adjoint operator $A$ has the $C$-symmetry, then
$C={J}e^Y$ (see (\ref{daad})) and the relation
$A{C}={C}A$ takes the form $A^*e^Y=e^YA$ (since $A^*J=JA$). This means
that the operator $A$ is similar to the self-adjoint operator $B$ given by
$$
B={e^{{Y}/2}}A{e^{{-Y}/2}}.
$$
\end{remark}

Let $S$ be a densely defined symmetric operator in the Hilbert space
$\mathfrak{H}$ with equal defect numbers. In what follows we suppose
that $S$ commutes with ${J}$, i.e,
\begin{equation}\label{neww31}
[S,J]=0.
\end{equation}

The condition (\ref{neww31}) naturally appears in the singular
perturbation theory when a singular (non-symmetric) potential
possesses a symmetry with respect ${J}$; see e.g. \cite{AK1}. 

\begin{definition}\label{daad43}
Let $S$ be a symmetric operator satisfying \eqref{neww31}. Denote by
$\{C_\alpha\}_{\alpha\in{D}}$ the collection of all possible
operators $C$ (parameterized by a set of indexes ${D}$) which
realize the property of $C$-symmetry for $S$ (in the sense of
Definition \ref{dad1}).
\end{definition}

According to (\ref{neww31}), ${J}\in\{C_\alpha\}_{\alpha\in{D}}$.

\begin{lemma}\label{ssss1}
Let $S$ be a symmetric operator satisfying \eqref{neww31}. Then the
following relation holds:
$$
C\in\{C_\alpha\}_{\alpha\in{D}}{\iff}C^*\in\{C_\alpha\}_{\alpha\in{D}}.
$$
\end{lemma}

\emph{Proof.} Let $C\in\{C_\alpha\}_{\alpha\in{D}}$. Then $S$
commutes with $C$ and $J$. Since $JC>0$, one obtains
$S(JC)=(JC)S=(C^*J)S$ and
\[
 S(C^*J)=S(JC)=(C^*J)S=C^*SJ,
\]
which yields $SC^*=C^*S$. Clearly, together with $C$, also $C^*$ has
the properties (i) and (ii) in Definition~\ref{dad1}. Thus, Lemma
\ref{ssss1} is proved. \rule{2mm}{2mm}

Denote by $\Sigma_{J}$ \emph{the set of all
${J}$-self-adjoint extensions of $S$}, i.e.,
$$
A\in\Sigma_{J}\ \iff \ {A\supset{S}} \quad \mbox{and} \quad
A^*J=JA.
$$

It follows from (\ref{neww31}) that
$A\in\Sigma_{J}\iff{A^*\in\Sigma_{J}}$.

\begin{definition}\label{dad2}\cite{KVS}
 We will say that an operator ${A}\in\Sigma_{J}$ belongs to the
set $\Sigma_{J}^{\textsf{st}}$ of stable ${C}$-symmetry if $A$ has
the property of ${C}$-symmetry realized by some operator
$C\in\{C_\alpha\}_{\alpha\in{D}}$.
\end{definition}

In other words, the condition $A\in\Sigma_{J}^{\textsf{st}}$ means
that the $C$-symmetry property of the operators $A$ and $S$ is
realized by the same operator $C$.

\begin{lemma}\label{neww20}
The following relation holds:
$A\in\Sigma_{J}^{\textsf{st}}\iff{A^*}\in\Sigma_{J}^{\textsf{st}}$.
\end{lemma}

\emph{Proof.} Let
$A\in\Sigma_{J}^{\textsf{st}}$. Then there exists an
operator $C\in\{C_\alpha\}_{\alpha\in{D}}$ commuting with $S$ and $A$.
By Lemma \ref{ssss1}, ${C}^{*}$ commutes with $S$.
On the other hand, ${C}^{*}$ commutes with $A^{*}$ by Proposition \ref{sese1}.
Therefore, $S$ and $A^{*}$ have the same ${C}^{*}$-symmetry. Lemma
\ref{neww20} is proved. \rule{2mm}{2mm}

\begin{remark}\label{fff1}
An arbitrary operator $C_\alpha$ from $\{C_\alpha\}_{\alpha\in{D}}$
determines a new definite inner product
\begin{equation}\label{dede1}
(\cdot,\cdot)_\alpha=[C_\alpha\cdot,\cdot]_{{J}}=({{J}}C_\alpha\cdot,\cdot)
\end{equation}
on the Krein space $(\mathfrak{H}, [\cdot, \cdot]_{{J}})$, which is
equivalent to the initial iner product $(\cdot,\cdot)=[{J}\cdot,
\cdot]_{{J}}$. The operator $S$ remains symmetric for any choice of
$(\cdot,\cdot)_\alpha$ and an arbitrary
$A\in\Sigma_{J}^{\textsf{st}}$ is, in fact, a self-adjoint extension
of $S$ with respect to a certain choice of $(\cdot,\cdot)_\alpha$.
\end{remark}

\subsection{The set $\Upsilon$.}
Denote by $\Upsilon$ the set of ${{J}}$-self-adjoint extensions
$A\in\Sigma_{J}$ which \emph{commute with every operator} from the
set $\{C_\alpha\}_{\alpha\in{D}}$, i.e.,
$$
A\in\Upsilon \quad \iff \quad  A\in\Sigma_{J} \quad
\mbox{and} \quad [A, C_\alpha]=0, \quad \forall\alpha\in{D}.
$$

Obviously,
$\Upsilon\subset\Sigma_{J}^{\textsf{st}}\subset\Sigma_{J}$.

\begin{theorem}\label{neww33}
The following statements are true:
\begin{itemize}
  \item[(i)] If $A\in\Upsilon$, then $A$ is a self-adjoint extension of
  $S$ for any choice of inner product $(\cdot,\cdot)_\alpha$.
  \item[(ii)] If $r$ is a real point of regular type of $S$, then the operator
  $A_r$ defined by (\ref{sas99}) belongs to $\Upsilon$.
  \item[(iii)] If $A\in\Upsilon$, then among all self-adjoint extensions
  $A'\supset{S}$ transversal to $A$ there exists at least one
  belonging to $\Upsilon$.
  \item[(iv)] If $S$ is nonnegative, then its Friedrichs and
  Krein-von Neumann extensions belong to $\Upsilon$.
\end{itemize}
\end{theorem}

\emph{Proof.} (i) Let $A\in\Upsilon$. Then $A$ is
${J}$-self-adjoint and
$C_\alpha{A}=AC_\alpha, \ \forall\alpha\in{D}$. This means that
$$
(Ax,y)_\alpha=[C_\alpha{A}x,y]_{{J}}=[C_\alpha{x},Ay]_{{J}}=(x,Ay)_\alpha, \quad \forall{x, y}\in\mathcal{D}(A).
$$
Therefore, $A$ is symmetric, and hence also self-adjoint, with
respect to $(\cdot,\cdot)_\alpha$.

(ii) It follows from Lemma \ref{ssss1} that $[S^*,C]=0$ for any
$C\in\{C_\alpha\}_{\alpha\in{D}}$. This implies that
\begin{equation}\label{neww63}
C(\ker(S^*-\mu{I}))=\ker(S^*-\mu{I}), \quad \forall\mu\in\mathbb{C}
\end{equation}
(since $C^2=I$). Using the definition (\ref{sas99}) of $A_r$ and
(\ref{neww63}) for $\mu=r$, we get $[A_r, C]=0$. Hence,
$A_r\in\Upsilon$.

(iii) According to the classical von Neumann formulas, an arbitrary
self-adjoint extension $A$ of $S$ is uniquely determined by a
unitary mapping $V:\ker(S^*+iI)\to\ker(S^*-iI)$: $x\in \cD(A)$ if
and only if $x_i=Vx_{-i}$, where $x_{\pm i}\in \ker(S^* \mp iI)$.
Let $A\in\Upsilon$. Then $A$ is self-adjoint (see (i)) and the
corresponding unitary mapping $V$ commutes with any
$C\in\{C_\alpha\}_{\alpha\in{D}}$ (due to (\ref{neww63})).
Considering the self-adjoint extension $A'\supset{S}$ determined by
the unitary mapping $V'=-V$ we establish (iii).

(iv) The Friedrichs extension $A_F$ and the Krein-von Neumann
extension $A_N$ of $S$ can be characterized as follows (see
\cite{AN} for the densely defined case and \cite{HSSW} for the
general case):

If $\{ f , f '\} \in S^{*}$, then $\{ f, f'\} \in A_{F}$ if and only
if
\begin{equation}\label{FR}
\inf \left \{ \| f - h \|^{2} + (f'- h' , f - h) \, : \, \{h,h'\}
\in S \right\} =0.
\end{equation}

If $\{ f , f '\} \in S^{*}$, then $\{ f, f'\} \in A_{N}$ if and only
if
\begin{equation}\label{VN}
\inf \left \{ \| f' - h' \|^{2} + (f'- h' , f - h) \, : \, \{h,h'\}
\in S \right\} =0.
\end{equation}

Since $[S^*,C]=[S,C]=0$ for any $C\in\{C_\alpha\}_{\alpha\in{D}}$,
formulas (\ref{FR}) and (\ref{VN}) imply that $A_F$ and $A_N$ are decomposed
with respect to the decomposition (\ref{d2}) (with subspaces $\mathfrak{L}_{\pm}$ determined by $C$):
$$
A_F=A_{F+}\dot{+}A_{F-}, \qquad A_N=A_{N+}\dot{+}A_{N-},
$$
where $A_{F\pm}$ and $A_{N\pm}$ are the Friedrichs extension and the
Krein-von Neumann extension of the symmetric operators
$S\upharpoonright{\mathfrak{L}_{\pm}}$ in the Hilbert spaces
$\mathfrak{L}_{\pm}$, respectively. These decompositions immediately
yield the relations $[A_F,C]=[A_N,C]=0$. Theorem \ref{neww33} is
proved. \rule{2mm}{2mm}

It may happen that the set $\Upsilon$ is empty, as the next example
shows.

\begin{example}\label{ex3}
Let $S_1$ and $S_2$ be symmetric operators with
defect numbers $<0,1>$ and $<1,0>$ acting in Hilbert spaces
$\mathfrak{H}_1$ and  $\mathfrak{H}_2$, respectively. Then
$S=S_1\oplus{S_2}$ is a symmetric operator in the Hilbert space
$\mathfrak{H}=\mathfrak{H}_1\oplus\mathfrak{H}_2$ with defect
numbers $<1,1>$. The operator $S$ commutes with the fundamental
symmetry ${J}=I\oplus{-I}$ in $\mathfrak{H}$.

Assume that $A\in\Upsilon$. Then $A$ is a self-adjoint extension of
$S$ (see Proposition \ref{neww30}) and $[A, {J}]=0$. Therefore,
$A\upharpoonright{\mathfrak{H}_1}$ is a self-adjoint extension of
$S_1$; however, this is impossible. Thus, $\Upsilon=\emptyset$.
\end{example}

\subsection{Description of $\Sigma_{J}^{\textsf{st}}$.}
The description of $\Sigma_{J}^{\textsf{st}}$ requires an appropriate boundary triplet
$(\mathcal{H}, \Gamma_0, \Gamma_1)$, in which the
images of $\{C_\alpha\}_{\alpha\in{D}}$ exist in the
parameter space $\mathcal{H}$.

\begin{lemma}\label{neww41}
For each $A\in\Upsilon$ there exists a boundary triplet
$(\mathcal{H}, \Gamma_0, \Gamma_1)$ for $S^*$ such that
$\mathcal{D}(A)=\ker\Gamma_0$ and the formulas
\begin{equation}\label{neww50}
{\mathcal C}_\alpha\Gamma_0f=\Gamma_0{C}_\alpha{f},  \quad {\mathcal
C}_\alpha\Gamma_1f=\Gamma_1{C}_\alpha{f}, \quad \forall{\alpha}\in{D}, \quad \forall{f}\in\mathcal{D}(S^*)
\end{equation}
correctly define the operator family $\{{\mathcal
C}_\alpha\}_{\alpha\in{D}}$ in $\mathcal{H}$. If $S$ is a simple
symmetric operator, then the correspondence
${{C}_\alpha}\rightarrow{\mathcal C}_\alpha$ established by
(\ref{neww50}) is injective.
\end{lemma}

\emph{Proof.} Let $A\in\Upsilon$. Then the corresponding unitary
mapping $V:\ker(S^*+iI)\to\ker(S^*-iI)$ (see the proof of Theorem
\ref{neww33}) commutes with the family
$\{C_\alpha\}_{\alpha\in{D}}$. Now introduce the boundary triplet
$(\mathcal{H}, \Gamma_0, \Gamma_1)$, where
$\mathcal{H}=\ker(S^*-iI)$ and
$$
\Gamma_0{x}=x_{i}-Vx_{-i}, \quad \Gamma_1{x}=ix_{i}+iVx_{-i}, \
x=u+x_{-i}+x_{i}\in\mathcal{D}(S^*),
$$
with $x_{\pm{i}}\in\ker(S^*{\mp}iI)$. Since $C_\alpha V=V C_\alpha$,
the restriction
$\mathcal{C}_\alpha=C_\alpha\upharpoonright{\ker(S^*-iI)}$ is an
operator in $\mathcal{H}$ and it satisfies the relations in
\eqref{neww50}. By construction, $\mathcal{D}(A)=\ker\Gamma_0$ (see
the proof of (iii) in Theorem~\ref{neww33}).

Let us assume that (\ref{neww50}) gives the same image $\mathcal{C}$
for two different $C$-symmetries $C_{\alpha_1}$ and $C_{\alpha_2}$
of $S$. Then
$(C_{\alpha_1}-C_{\alpha_2})\mathcal{D}(S^*)\subset\mathcal{D}(S)$.
Combining this relation with  (\ref{neww63}), one concludes that
$C_{\alpha_1}f_\mu=C_{\alpha_2}f_\mu$ for every
$f_\mu\in\ker(S^*-\mu{I})$ and
$\mu\in\mathbb{C}\setminus\mathbb{R}$. Since the symmetric operator
$S$ is simple, this implies that $C_{\alpha_1}={C_{\alpha_2}}$.
Lemma \ref{neww41} is proved. \rule{2mm}{2mm}

\begin{remark}\label{remark1}
If a boundary triplet $(\mathcal{H}, \Gamma_0, \Gamma_1)$ of $S^*$
satisfies (\ref{neww50}), then the associated transversal extensions
$A_0$ and $A_1$ defined by (\ref{nnn1}) belong to $\Upsilon$. This
means that \emph{boundary triplets $(\mathcal{H}, \Gamma_0,
\Gamma_1)$ with property (\ref{neww50}) exist if and only if the set
$\Upsilon$ is non-empty}.
\end{remark}

The existence of a boundary triplet $(\mathcal{H}, \Gamma_0,
\Gamma_1)$ for $S^*$, which satisfies the properties in
\eqref{neww50}, will guarantee a couple of useful properties for the
operators ${\mathcal C}_\alpha$ and also some important relations
between the extensions of $S$ and the parameters corresponding to
them in $\mathcal{H}$.

\begin{lemma}\label{newlem}
Let $(\mathcal{H}, \Gamma_0, \Gamma_1)$ be a boundary triplet for
$S^*$ with the properties (\ref{neww50}). Then the associated
$\gamma$-field $\gamma(\cdot)$ and Weyl function $M(\cdot)$ satisfy
the following relations
\begin{equation}\label{new01}
 \gamma(\lambda) \mathcal{C}_\alpha = C_\alpha \gamma(\lambda)
 \quad \text{ and } \quad [\mathcal{C}_\alpha, M(\lambda)]=0
\end{equation}
for all $\lambda\in\mathbb{C}_+\cup\mathbb{C}_-$ and $\alpha\in D$.
\end{lemma}
\emph{Proof.} Since
$\mathcal{C}_{\alpha}(\ker(S^*-\lambda))=\ker(S^*-\lambda)$, see
\eqref{neww63}, both of the identities in \eqref{new01} follow
easily from (\ref{neww50}) by applying the formula
$\gamma(\lambda)=(\Gamma_0\upharpoonright \ker(S^*-\lambda))^{-1}$
and the definition of $M(\lambda)$ in \eqref{neww65}.
\rule{2mm}{2mm}

The next lemma concern the class of operators $\{{\mathcal
C}_\alpha\}_{\alpha\in D}$ appearing in Lemma~\ref{neww41}. In the
case $C_\alpha=J$, we will use the special notation ${\mathcal J}$
for the corresponding intertwining operator ${\mathcal C}_\alpha$ in
(\ref{neww50}), i.e., ${\mathcal J}\Gamma_0=\Gamma_0{J}$ and
${\mathcal J}\Gamma_1=\Gamma_1{J}$.

\begin{lemma}\label{neww67}
Let $(\mathcal{H}, \Gamma_0, \Gamma_1)$ be a boundary triplet of
$S^*$ with the properties (\ref{neww50}). Then the intertwining
operators ${\mathcal C}_\alpha$ satisfy the relations
${{\mathcal{C}}}_\alpha^2=I,$ \  $\mathcal{JC}_\alpha>0$ and the
operator ${\mathcal J}$ is a fundamental symmetry in the auxiliary
Hilbert space $\mathcal{H}$.
\end{lemma}
\emph{Proof.} The identity ${{\mathcal{C}}}_\alpha^2=I$ immediately
follows from (\ref{neww50}) and the corresponding identity
${{{C}}}_\alpha^2=I$. Furthermore, $[S^*, {J}]=0$ due to
(\ref{neww31}). Taking this relation into account and considering
\eqref{new2} (or \eqref{new22}, see \eqref{new21}) with ${J}x$ and
${J}y$ instead of $x$ and $y$ it follows that ${\mathcal J}$ is
unitary in $\mathcal{H}$. This together with the identity ${\mathcal
J}^2=I$ leads to the fact that ${\mathcal J}$ is a self-adjoint
operator. Hence, ${\mathcal J}$ is a fundamental symmetry in
$\mathcal{H}$.

Since $[S^*,C_\alpha]=0$ for every $\alpha\in{D}$, one can rewrite
(\ref{new22}) by substituting ${J}C_\alpha{x}$ instead of $x$ as
follows:
\begin{equation}\label{neww64}
(S^*x,y)_\alpha-(x,S^*y)_\alpha=i[({\mathcal
J}\mathcal{C}_\alpha\Omega_{+}x,\Omega_{+}y)_{\mathcal H}-({\mathcal
J}\mathcal{C}_\alpha\Omega_{-}x,\Omega_{-}y)_{\mathcal H}],
\end{equation}
where $(\cdot,\cdot)_\alpha=({{J}}C_\alpha\cdot, \cdot)$. Now by
putting $x=y={f}_\mu\in\ker(S^*-\mu{I})$, $\mu\in\mathbb{C}_+$, in
(\ref{neww64}) and recalling the definition (\ref{neww65}) of the
characteristic function $\Theta({\mu})$ one obtains
\begin{equation}\label{neww70}
2(\textsf{Im} \ \mu)({f}_\mu, {f}_\mu)_{\alpha}=({\mathcal
J}\mathcal{C}_\alpha{h},h)_{\mathcal H}-({\mathcal
J}\mathcal{C}_\alpha\Theta({\mu})h, \Theta({\mu})h)_{\mathcal H},
\end{equation}
where $h=\Omega_{+}{f}_\mu$ is an arbitrary element of ${\mathcal
H}$ (since $\Omega_{+}$ maps $\ker(S^*-\mu{I})$ onto ${\mathcal
H}$). Due to (\ref{neww65})  and (\ref{neww50}), $[\Theta({\mu}),
\mathcal{C}_\alpha]=0 \ (\forall\mu\in\mathbb{C}_\pm)$. Hence,
(\ref{neww70}) implies that
 $$
({\mathcal
J}\mathcal{C}_\alpha(I-\Theta^*({\mu})\Theta({\mu})){h},h)_{\mathcal
H}=({\mathcal J}\mathcal{C}_\alpha{Fh},Fh)_{\mathcal H}>0,
$$
where $F=(I-\Theta^*({\mu})\Theta({\mu}))^{1/2}$ is an invertible
operator ($0\in\rho(F)$) in ${\mathcal H}$ (due to
$\|\Theta({\mu})\|=\|\Theta^*({\mu})\|<1$) such that
$[F,\mathcal{C}_\alpha]=0$. Therefore ${\mathcal
J}\mathcal{C}_\alpha>0$. Lemma \ref{neww67} is proved.
\rule{2mm}{2mm}

\begin{theorem}\label{neww80}
Let $({\mathcal H}, \Gamma_0, \Gamma_1)$ be a boundary triplet of
$S^*$ with the properties (\ref{neww50}). Then an arbitrary
$A\in\Sigma_{J}^{\textsf{st}}$ admits the presentation
\begin{equation}\label{as1}
A=S^*\upharpoonright{\{f\in\mathcal{D}(S^*) \ | \
\mathcal{K}\Omega_{+}f=\Omega_{-}f\}},
\end{equation}
where $\Omega_\pm$ are defined by (\ref{new21}) and $\mathcal{K}$ is a stable ${\mathcal J}$-unitary operator in
the Krein space $({\mathcal H}, [\cdot,\cdot]_{{\mathcal J}})$.
\end{theorem}

\emph{Proof.} It follows from (\ref{daad}) and Lemma \ref{neww67}
that ${\mathcal J}\mathcal{C}_\alpha=e^{{\mathcal Y}_\alpha}$, where
${\mathcal Y}_\alpha$ is a bounded self-adjoint operator in
$\mathcal{H}$ and $\{{\mathcal J}, {\mathcal Y}_\alpha\}=0$.

Since the operator $S^*$ remains to be adjoint for the symmetric
operator $S$ with respect to $(\cdot,\cdot)_\alpha$, formulas
\eqref{new21}, (\ref{new22}), and (\ref{neww64}) imply that
$({\mathcal H}, e^{{\mathcal Y}_\alpha/2}\Gamma_0, e^{{\mathcal
Y}_\alpha/2}\Gamma_1)$ is a boundary triplet for $S^*$ acting in the
Hilbert space ${\mathfrak H}$ with the inner product $(\cdot,
\cdot)_{\alpha}$.

Let $A\in\Sigma_{J}^{\textsf{st}}$. Then $[A, C_\alpha]=0$ for a
certain choice of ${\alpha\in{D}}$ and $A$ is a self-adjoint
extension of $S$ with respect to $(\cdot,\cdot)_\alpha$. Therefore,
$A$ has the presentation (see \cite{Gor}):
$$
A=S^*\upharpoonright{\{f\in\mathcal{D}(S^*) \ | \
\mathcal{W}e^{\mathcal{Y}_\alpha/2}\Omega_{+}f=e^{\mathcal{Y}_\alpha/2}\Omega_{-}f\}},
$$
where $\mathcal{W}$ is a unitary operator in the Hilbert space
$\mathcal{H}$. The previous description of $A$ leads to (\ref{as1})
with
\begin{equation}\label{neww909b}
\mathcal{K}=e^{-\mathcal{Y}_\alpha/2}\mathcal{W}e^{\mathcal{Y}_\alpha/2}.
\end{equation}
Since the adjoint\footnote{$A^*$ denotes the adjoint with respect to
the initial scalar product $(\cdot,\cdot)$ of $\mathfrak{H}$}
operator $A^*$ is determined by (\ref{as1}) with ${{\mathcal
K}^*}^{-1}$ (cf. \cite{Gor}), the relation $A^*{J}={J}A$ means that
${{\mathcal K}^*}^{-1}{\mathcal J}={\mathcal J}{\mathcal K}$.
Therefore, ${\mathcal K}$ is a ${\mathcal J}$-unitary operator in
the Krein space $({\mathcal H}, [\cdot,\cdot]_{{\mathcal J}})$.

It follows from (\ref{neww909b}) that $\|{\mathcal
K}^n\|<\mbox{const}, \ \forall{n}\in\mathbb{Z}$. Hence, ${\mathcal
K}$ is a stable ${\mathcal J}$-unitary operator \cite{AZ}. Theorem
\ref{neww80} is proved. \rule{2mm}{2mm}

\begin{corollary}\label{ess1}
The formula (\ref{as1}) determines an operator $A\in\Sigma_{J}^{\textsf{st}}$ if and only if
the corresponding ${\mathcal J}$-unitary operator ${\mathcal K}$ has the property of
${C}$-symmetry realized by an operator $\mathcal{C}_\alpha$ from $\{{\mathcal C}_\alpha\}_{\alpha\in{D}}$.
\end{corollary}
\emph{Proof.} If $A\in\Sigma_{J}^{\textsf{st}}$ then
$[A, C_\alpha]=0$ for a certain choice of ${\alpha\in{D}}$. This means that $C_\alpha : \mathcal{D}(A)\to\mathcal{D}(A)$.
By (\ref{new21}) and (\ref{neww50}), ${\mathcal C}_\alpha\Omega_\pm=\Omega_\pm{C_\alpha}$. Combining this with
(\ref{as1}) and taking into account that ${\mathcal C}_\alpha^2=I$, we obtain $[{\mathcal K}, {\mathcal C}_\alpha]=0$.

Conversely, if $A$ is determined by (\ref{as1}) and $[{\mathcal K},
{\mathcal C}_\alpha]=0$ for a certain operator ${\mathcal C}_\alpha$
from $\{{\mathcal C}_\alpha\}_{\alpha\in{D}}$, then for its
`preimage' $C_\alpha$ (see (\ref{neww50})) the relation $C_\alpha :
\mathcal{D}(A)\to\mathcal{D}(A)$ holds. This means that $[A,
C_\alpha]=0$ (since $[S^*, C_\alpha]=0$). Hence,
$A\in\Sigma_{J}^{\textsf{st}}$. \rule{2mm}{2mm}

\begin{remark}\label{wawe1}
In general, the formula (\ref{as1}) establishes a bijection between
the elements of $\Sigma_{J}^{\textsf{st}}$ and some subset of the
set of stable ${\mathcal J}$-unitary operators $\mathcal{K}$ in
$({\mathcal H}, [\cdot,\cdot]_{{\mathcal J}})$. This subset is
uniquely determined by the additional assumption that ${\mathcal K}$
has the property of ${C}$-symmetry realized by an operator from the
image ${\{{\mathcal C}_\alpha\}_{\alpha\in{D}}}$ of the set
${\{{C}_\alpha\}_{\alpha\in{D}}}$ (see (\ref{neww50})). However, it
is not easy to apply such kind of definition. A more appropriate
external description for ${\{{\mathcal C}_\alpha\}_{\alpha\in{D}}}$
is established in the next subsection by means of reproducing kernel
Hilbert space models associated with Nevanlinna functions.
\end{remark}

\subsection{Reproducing kernel Hilbert space models.}
Let $M(\cdot)$ be a Weyl function of the symmetric operator $S$
associated with the boundary triplet $(\mathcal{H}, \Gamma_0,
\Gamma_1)$. The corresponding \emph{Nevanlinna kernel}
$\textsf{N}_M(\xi,\mu)$ on
$(\mathbb{C}_+\cup\mathbb{C}_-)\times(\mathbb{C}_+\cup\mathbb{C}_-)$
is defined by
\begin{equation}\label{neww90}
\textsf{N}_M(\xi,\mu):=
\frac{M(\mu)-M(\overline{\xi})}{\mu-\overline{\xi}},
\quad \mu,\xi\in\mathbb{C}_+\cup\mathbb{C}_-, \quad \xi\not=\overline{\mu}.
\end{equation}
The kernel $\textsf{N}_M(\xi,\mu)$ is Hermitian, holomorphic, and nonnegative.
The corresponding reproducing kernel Hilbert space will be denoted by $\mathfrak{H}_M$.
The space $\mathfrak{H}_M$ consists of $\mathcal{H}$-valued holomorphic vector functions on
$\mathbb{C}_+\cup\mathbb{C}_-$ obtained as the closed linear span of functions $\mu\to\textsf{N}_M(\xi,\mu)f,$ \
$\xi\in\mathbb{C}_+\cup\mathbb{C}_-, \ f\in\mathcal{H}$, which is provided with the scalar product determined by
$$
\prec\textsf{N}_M(\xi,\cdot)f, \textsf{N}_M(\lambda,\cdot)g\succ:=
 (\textsf{N}_M(\xi,\lambda)f,g)_{\mathcal{H}}, \quad f,g\in\mathcal{H},
 \quad \xi, \lambda\in\mathbb{C}_\pm.
$$
The functions $\phi(\cdot)\in\mathfrak{H}_M$ satisfy the reproducing kernel property
\begin{equation}\label{neww92}
\prec\phi(\cdot),
\textsf{N}_M(\lambda,\cdot)g\succ=(\phi(\lambda),g)_{\mathcal{H}},
\quad g\in\mathcal{H}, \quad \lambda\in\mathbb{C}_\pm.
\end{equation}

The reproducing kernel Hilbert space $\mathfrak{H}_M$ gives rise to
a useful model representation of the symmetric operator $S$ and the
associated boundary mappings. The next statement contains a lot of
relevant results; see \cite{BHS} for a detailed proof and further
details.
\begin{proposition}\label{neww95}
Let $M(\cdot)$ be a Weyl function of a simple symmetric operator $S$. Then:
\begin{itemize}
  \item[(i)] the linear relation $S_M=\{\{\phi,\psi\}\in\mathfrak{H}_M^2 : \psi(\lambda)=\lambda\phi(\lambda)\}$ is a symmetric
  operator in $\mathfrak{H}_M$ which is unitarily equivalent to $S$;
  \item[(ii)] the linear relation
  $$
  \mathcal{T}=\{\{\phi,\psi\}\in\mathfrak{H}_M^2 : \psi(\lambda)-\lambda\phi(\lambda)=c_1+M(\lambda)c_2, \ c_1,c_2\in\mathcal{H}\}
  $$
  determines the adjoint $S_M^*$ of $S_M$ in $\mathfrak{H}_M$;
  \item[(iii)] The operators
  $$
  \Gamma_0^M\{\phi,\psi\}=c_2, \qquad \Gamma_1^M\{\phi,\psi\}=-c_1, \qquad \{\phi,\psi\}\in\mathcal{T}
  $$
  form a boundary triplet $(\mathcal{H}, \Gamma_0^M, \Gamma_1^M)$ for $S_M^*$;
  \item[(iv)] The Weyl function of $S_M$ associated with $(\mathcal{H}, \Gamma_0^M, \Gamma_1^M)$ coincides with
  $M(\cdot)$.
\end{itemize}
\end{proposition}

Proposition \ref{neww95} allows one to establish an `external' description of
the set ${\{{\mathcal C}_\alpha\}_{\alpha\in{D}}}$ (see Remark \ref{wawe1}).

\begin{theorem}\label{ttt6}
Let $S$ be a simple symmetric operator, let $({\mathcal H},
\Gamma_0, \Gamma_1)$ be a boundary triplet for $S^*$ with the
properties (\ref{neww50}), and let $M(\cdot)$ be the associated Weyl
function. Then a bounded operator $\mathcal{C}$ in $\mathcal{H}$
belongs to the set ${\{{\mathcal C}_\alpha\}_{\alpha\in{D}}}$ if and
only if
\begin{equation}\label{wawe2}
(i)  \ \mathcal{C}^2=I; \quad (ii)  \ \mathcal{JC}>0; \quad (iii)  \
[\mathcal{C}, M(\lambda)]=0, \ \forall\lambda\in\mathbb{C}_\pm.
\end{equation}
\end{theorem}

\emph{Proof.} If $\mathcal{C}\in{\{{\mathcal
C}_\alpha\}_{\alpha\in{D}}}$, then (\ref{wawe2}) holds by
Lemma~\ref{newlem} and Lemma~\ref{neww67}.

Now the converse will be proved. If $\mathcal{C}$ satisfies
(\ref{wawe2}), then its adjoint $\mathcal{C}^*$ also satisfies
(\ref{wawe2}). This means that the operators $\textsf{C}$ and
$\textsf{C}'$,
$$
\textsf{C}[\textsf{N}_M(\lambda,\cdot)f]:=\textsf{N}_M(\lambda,\cdot)\mathcal{C}f,
\qquad
\textsf{C}'[\textsf{N}_M(\lambda,\cdot)f]:=\textsf{N}_M(\lambda,\cdot)\mathcal{C}^*f,
$$
are correctly defined on the linear span of functions
$\{\textsf{N}_M(\lambda,\cdot)f\}$ with
${\lambda\in\mathbb{C}_+\cup\mathbb{C}_-, f\in\mathcal{H}}$. It
follows from \eqref{neww92} that
\begin{eqnarray*}
 \|\textsf{N}_M(\lambda,\cdot)\mathcal{C}f\|^2_{\sH_M}
 &=&(\textsf{N}_M(\lambda,\lambda)\mathcal{C}f,\mathcal{C}f)_\cH \vspace{4mm} \\
  &\le&\|{C}\|^2_\cH\|(\textsf{N}_M(\lambda,\lambda)f,f)_\cH \vspace{4mm} \\
  &=&\|{C}\|^2_\cH\|\textsf{N}_M(\lambda,\cdot)\mathcal{C}f\|^2_{\sH_M}.
\end{eqnarray*}
Thus $\textsf{C}$ and, similarly, $\textsf{C}'$ is continuous.
Hence, $\textsf{C}$ and $\textsf{C}'$ can be extended by continuity
onto the whole space $\mathfrak{H}_M$. Using (\ref{neww92}) twice
gives
$$
\prec\phi(\cdot), \textsf{N}_M(\lambda,\cdot)\mathcal{C}^*g\succ=(\phi(\lambda),\mathcal{C}^*g)_{\mathcal{H}}=(\mathcal{C}\phi(\lambda),g)_{\mathcal{H}}
$$
and
\begin{eqnarray*}
 \prec\phi(\cdot),  \textsf{N}_M(\lambda,  \cdot)\mathcal{C}^*g\succ
 =\prec\phi(\cdot), \textsf{C}'[\textsf{N}_M(\lambda,\cdot)g]\succ   \vspace{4mm} \\
 = \prec{\textsf{C}'}^*[\phi(\cdot)], \textsf{N}_M(\lambda,\cdot)g\succ
 &=(({\textsf{C}'}^*[\phi(\cdot)])(\lambda),g)_{\mathcal{H}}.
\end{eqnarray*}

Comparing the righthand sides, one obtains
$({\textsf{C}'}^*[\phi(\cdot)])(\lambda)=\mathcal{C}\phi(\lambda)$
for all $\lambda\in\mathbb{C}_+\cup\mathbb{C}_-.$ This means that
${\textsf{C}'}^*=\textsf{C}$ and that the action of $\textsf{C}$ on
an arbitrary $\mathcal{H}$-valued function
$\phi(\cdot)\in\mathfrak{H}_M$ is realized via the action of
$\mathcal{C}$ on the vectors
 $\phi(\lambda)\in{\mathcal H}$, i.e.
\begin{equation}\label{uman2}
 ({\textsf{C}}[\phi(\cdot)])(\lambda)\equiv\mathcal{C}\phi(\lambda),
 \quad  \forall\lambda\in\mathbb{C}_\pm.
 \end{equation}
Therefore, $\textsf{C}^2=I$ (since $\mathcal{C}^2=I$) and it is
clear from (i) in Proposition~\ref{neww95} that $\textsf{C}$
commutes with $S_M$.

Repeating the arguments above for the case where
$\mathcal{C}=\mathcal{J}$, one obtains a fundamental symmetry
$\textsf{J}$ in $\mathfrak{H}_M$ defined by the formula
\begin{equation}\label{wewe1}
(\textsf{J}[\phi(\cdot)])(\lambda)\equiv\mathcal{J}\phi(\lambda).
\end{equation}
It is easy to see that $\textsf{JC}>0$ in $\mathfrak{H}_M$, since
$\mathcal{JC}>0$ in $\mathcal{H}$. Thus, starting with an operator
$\mathcal{C}$ satisfying (\ref{wawe2}), one can construct the
operator ${\textsf{C}}$, which realizes the property of
${C}$-symmetry for the symmetric operator $S_M$ in the Krein space
$(\mathfrak{H}_M, [\cdot, \cdot]_{\textsf{J}})$.

By Proposition \ref{neww95}, $S_M$ and $S$ have the same Weyl
function $M(\cdot)$ associated with the boundary triplets
$({\mathcal H}, \Gamma_0^M, \Gamma_1^M)$ and $({\mathcal H},
\Gamma_0, \Gamma_1)$, respectively. Therefore, there exists a
unitary mapping $U :
\mathfrak{H}_M\stackrel{\textrm{onto}}{\to}\mathfrak{H}$  such that
$S_M=U^{-1}SU$ and $\Gamma_j^M=\Gamma_jU$, $j=0,1$; see
\cite[Theorem~3.9]{DHMS}.

Let us show that $U$ can be chosen in such a way that
\begin{equation}\label{wawe3}
 \textsf{J}=U^{-1}{J}U,
\end{equation}
where $\textsf{J}$ is defined by (\ref{wewe1}). Indeed, it follows
from (\ref{neww31}) and (\ref{neww50}) that
$M(\cdot)=M_+(\cdot)\oplus{M_-(\cdot)}$, where the decomposition is
with respect to the fundamental decomposition
$\mathcal{H}=\mathcal{H}_+\oplus\mathcal{H}_-$ of the Krein space
$(\mathcal{H}, [\cdot,\cdot]_{\mathcal{J}})$. Furthermore,
$M_\pm(\cdot)$ are the Weyl functions of the symmetric operators
$S_\pm=S\upharpoonright{\mathfrak{H}_\pm}$ acting in subspaces
$\mathfrak{H}_\pm$ of the fundamental decomposition (\ref{d1}) of
the Krein space $(\mathfrak{H}, [\cdot,\cdot]_{{J}})$.

Let $\mathfrak{H}_{M_{\pm}}$ be the reproducing kernel Hilbert
spaces constructed by $M_\pm(\cdot)$. In view of (\ref{wewe1})
$\mathfrak{H}_{M}=\mathfrak{H}_{M_{+}}\oplus\mathfrak{H}_{M_{-}}$ is
the fundamental decomposition of the Krein space $(\mathfrak{H}_M,
[\cdot,\cdot]_{\textsf{J}})$ and one has
$S_M=S_{M_+}\oplus{S}_{M_-}$ with respect to this decomposition. The
pairs of operators $S_{M_+}$, $S_+$ and $S_{M_-}$, $S_-$ have the
Weyl functions $M_+(\cdot)$ and $M_-(\cdot)$, respectively. These
functions $M_{\pm}(\cdot)$ are associated with the boundary triplets
$({\mathcal H}_{\pm}, \Gamma_0^{M_\pm}, \Gamma_1^{M_\pm})$ and
$({\mathcal H}_\pm, \Gamma_0^\pm, \Gamma_1^\pm)$ of $S_{M_\pm}^*$
and $S_\pm^*$, respectively. Here $\Gamma_j^{M_\pm}$ are defined
according to the statement (iii) of Proposition \ref{neww95} and
$\Gamma_j^\pm$ are the restrictions of $\Gamma_j$ onto
$\mathcal{D}(S_{\pm}^*)$. Without loss of generality, one can choose
unitary mappings $U_\pm :
\mathfrak{H}_{M_\pm}\stackrel{\textrm{onto}}{\to}\mathfrak{H}_{\pm}$
such that $S_{M_\pm}=U_\pm^{-1}S_{\pm}U_\pm$ and
$\Gamma_j^{M_\pm}=\Gamma_j^\pm{U_\pm}$. But, then the operator
$U=U_+\oplus{U_-}$ satisfies (\ref{wawe3}).

It follows from (\ref{wawe3}) that the set
$\mathcal{U}=\{U^{-1}C_\alpha{U}\}_{\alpha\in{D}}$ contains all
possible $C$-symmetries of $S_M$ in the Krein space
$(\mathfrak{H}_M, [\cdot, \cdot]_{\textsf{J}})$. Therefore, the
operator $\textsf{C}$ defined by (\ref{uman2}) belongs to
$\mathcal{U}$ and $\textsf{C}=U^{-1}C_\alpha{U}$ for a certain
choice of $\alpha\in{D}$. In that case, taking (\ref{neww50}) into
account, one obtains
$$
\Gamma_j^M\textsf{C}=\Gamma_j^MU^{-1}C_\alpha{U}=\Gamma_jC_\alpha{U}={\mathcal C}_{\alpha}\Gamma_j{U}={\mathcal C}_{\alpha}\Gamma_j^M \quad
(j=0,1).
$$

On the other hand, in view of (\ref{uman2}) and the statements (ii),
(iii) in Proposition \ref{neww95} one has
$\Gamma_j^M\textsf{C}={\mathcal C}\Gamma_j^M$. Consequently,
${\mathcal C}={\mathcal C}_\alpha$. Theorem \ref{ttt6} is proved.
\rule{2mm}{2mm}

\subsection{Resolvent formula for $J$-self-adjoint extensions with stable $C$-symmetry.}

Combining Theorem \ref{ttt6} with Corollary \ref{ess1} we immediately
obtain the following complete description of
$\Sigma_{J}^{\textsf{st}}$:

\begin{theorem}\label{neww80b}
Let $({\mathcal H}, \Gamma_0, \Gamma_1)$ be a boundary triplet of
$S^*$ with properties (\ref{neww50}) and let $M(\cdot)$ be the Weyl
function of $S$. Then $A\in\Sigma_{J}^{\textsf{st}}$ if and only if
$A$ is defined by (\ref{as1}) and the corresponding ${\mathcal
J}$-unitary operator $\mathcal{K}$ has the $\mathcal{C}$-symmetry in
$({\mathcal H}, [\cdot,\cdot]_{{\mathcal J}})$ such that
$[\mathcal{C}, M(\cdot)]=0$.
\end{theorem}

Another characterization for $A$ to belong to the class
$\Sigma_{J}^{\textsf{st}}$ can be obtained by describing the
resolvents of $A\in\Sigma_{J}^{\textsf{st}}$. Recall from
Remark~\ref{fff1} that if $A\in\Sigma_{J}^{\textsf{st}}$, then $A$
is selfadjoint in the Hilbert space $(\sH,(\cdot,\cdot)_\alpha)$,
where the inner product is defined by \eqref{dede1}. Therefore, the
resolvent set of the $J$-self-adjoint operator
$A\in\Sigma_{J}^{\textsf{st}}$ is automatically nonempty, since
$\dC_\pm\subset\rho(A)$. To establish such a characterization, the
following definition is needed.

\begin{definition}\label{commdef}
Let $R$ be a (closed linear) relation in a Hilbert space $\cH$ and
let $C$ be a bounded operator in $\cH$. Then $C$ is said to commute
with $R$, if the following formula holds:
\begin{equation}\label{comm1}
  R=\{\{Cf,Cf'\}:\, \{f,f'\}\in R\}.
\end{equation}
In this case we write shortly $[R,C]=0$.
\end{definition}

Observe, that if $R$ is an operator then $\{f,f'\}\in R$ means that
$f'=Rf$, and thus \eqref{comm1} can be rewritten as $RCf=CRf$ for
all $f\in\dom R$, i.e., Definition~\ref{commdef} reduces to the
usual definition of commutativity, when $R$ is an operator. Indeed,
it is straightforward to check that the condition \eqref{comm1} is
equivalent to $RC=CR$, where the products are to be understood in
the sense of relations.

In the next statement $\gamma(\cdot)$ stands for the $\gamma$-field
corresponding to the boundary triplet $({\mathcal H}, \Gamma_0,
\Gamma_1)$ in Theorem \ref{neww80b} and
$A_0=S^*\upharpoonright\ker\Gamma_0$.

\begin{theorem}\label{absd1}
Let the assumptions be as in Theorem~\ref{neww80b}. Then
$A\in\Sigma_{J}^{\textsf{st}}$ if and only if
\begin{equation}\label{bbb1}
(A-\mu{I})^{-1}=(A_0-\mu{I})^{-1}-\gamma(\mu)(M(\mu)-\mathcal{R})^{-1}\gamma^*(\overline{\mu}),
 \quad \mu\in \rho(A)\cap\rho(A_0),
\end{equation}
where $\mathcal{R}$ is a ${\mathcal J}$-self-adjoint
relation\footnote{We refer to \cite{Der95} for the basic definitions
of linear relation theory in the Krein space setting.}, which has
the $\mathcal{C}$-symmetry in $({\mathcal H},
[\cdot,\cdot]_{{\mathcal J}})$ such that $[\mathcal{C},M(\cdot)]=0$.

Furthermore, $A$ is disjoint with $A_0$ (i.e.,
$\mathcal{D}(A)\cap\mathcal{D}(A_0)=\mathcal{D}(S)$) if and only if
$\mathcal{R}$ is an operator with the indicated
$\mathcal{C}$-symmetry, and $A$ is transversal with $A_0$ (i.e.,
$\mathcal{D}(A)\dot{+}\mathcal{D}(A_0)=\mathcal{D}(S^*)$) if and
only if $\mathcal{R}$ is a bounded operator with the indicated
$\mathcal{C}$-symmetry.
\end{theorem}
\emph{Proof.} First assume that $A\in\Sigma_{J}^{\textsf{st}}$. Then
there exists $C_\alpha\in\{C_\alpha\}_{\alpha\in{D}}$ such that
$[A,C_\alpha]=0$. Moreover, $[A_0,C_\alpha]=0$ (since
$A_0\in\Upsilon$, see Remark \ref{remark1}). This mean that $A$ and
$A_0$ are self-adjoint extensions of the symmetric operator $S$ with
respect to the scalar product $(\cdot,\cdot)_\alpha$ (Remark
\ref{fff1}) and hence, in particular, $\dC_\pm\subset
\rho(A)\cap\rho(A_0)$.

Now, rewrite (\ref{as1}) as follows
\begin{equation}\label{as1bbb}
 A=S^*\upharpoonright{\{f\in\mathcal{D}(S^*) :\,
 i(I+\mathcal{K})\Gamma_0f =(I-\mathcal{K})\Gamma_1f \}}.
\end{equation}
Since $\mathcal{K}$ is $J$-unitary, this means that $A$ corresponds
to the ${\mathcal J}$-self-adjoint relation
$\mathcal{R}=i(I+\mathcal{K})(I-\mathcal{K})^{-1}$ in $\cH$, i.e.,
$\{\Gamma_0,\Gamma_1\}\cD(A)=\mathcal{R}$. By Theorem~\ref{neww80b}
$\mathcal{K}$ has the $\mathcal{C}$-symmetry in $({\mathcal H},
[\cdot,\cdot]_{{\mathcal J}})$ realized by an operator
$\mathcal{C}\in\{\mathcal{C}_\alpha\}_{\alpha\in D}$ such that
$[\mathcal{C},M(\cdot)]=0$. Since
$\{h,k\}=\{\Gamma_0f,\Gamma_1f\}\in \mathcal{R}$ if and only if
$i(I+\mathcal{K})\Gamma_0f =(I-\mathcal{K})\Gamma_1f$, see
\eqref{as1bbb}, and $[\mathcal{C},\mathcal{K}]=0$, it is clear that
\eqref{comm1} is satisfied, so that $[\mathcal{C},\mathcal{R}]=0$.
This means that $\mathcal{R}$ has the $\mathcal{C}$-symmetry in
$({\mathcal H}, [\cdot,\cdot]_{{\mathcal J}})$.

Finally, since $\mathcal{R}$ corresponds to $A$
($\{\Gamma_0,\Gamma_1\}\cD(A)=\mathcal{R}$) and
$\dC_\pm\subset\rho(A)$, \cite[Proposition~2.1]{DM} shows that
$0\in\rho(M(\mu)-\mathcal{R})$ for all $\mu\in \dC_\pm$, and
moreover, the resolvent formula \eqref{bbb1} is obtained from
\cite[Proposition~2.2]{DM} (notice that these two Propositions in
\cite{DM} are formulated for an arbitrary closed linear relation
$\mathcal{R}$ in $\cH$).

To prove the converse statement assume that $A$ is given by
\eqref{bbb1} for some $J$-selfadjoint relation $\mathcal{R}$, which
has the $\mathcal{C}$-symmetry in $({\mathcal H},
[\cdot,\cdot]_{{\mathcal J}})$ such that $[\mathcal{C},M(\cdot)]=0$.
Then $[\mathcal{C},\mathcal{R}]=0$ and this is equivalent to
$[\mathcal{C},(\mathcal{R}-M(\mu))^{-1}]=0$, $\mu\in\dC_\pm$, since
$[\mathcal{C},M(\cdot)]=0$. To see this, observe that $\{f,f'\}\in
\mathcal{R}$ is equivalent to
\[
 \{f'-M(\mu)f,f\}\in (\mathcal{R}-M(\mu))^{-1}
\]
and that $\{\mathcal{C}f,\mathcal{C}f'\}\in \mathcal{R}$ is
equivalent to
\[
 \{\mathcal{C}f'-M(\mu)\mathcal{C}f,\mathcal{C}f\}
 =\{\mathcal{C}(f'-M(\mu)f),\mathcal{C}f\}
 \in (\mathcal{R}-M(\mu))^{-1},
\]
since $[\mathcal{C},M(\cdot)]=0$. Here $(\mathcal{R}-M(\mu))^{-1}$
is an operator for $\mu\in\dC_\pm$, and hence we conclude that
$[\mathcal{C},\mathcal{R}]=0$ if and only if $\mathcal{C}$ commutes
with $(\mathcal{R}-M(\mu))^{-1}$. Now it follows from
Lemma~\ref{newlem} that the 'preimage' $C_\alpha$ of
$\mathcal{C}=\mathcal{C}_\alpha$, $\alpha\in D$, (see
\eqref{neww50}) commutes with both of the summands in the righthand
side of \eqref{bbb1}. Thus, $[C_\alpha,(A-\mu I)^{-1}]=0$ and this
is equivalent to $[C_\alpha,A]=0$. Furthermore, since $J\in
\{C_\alpha\}_{\alpha\in{D}}$ it is clear from \eqref{bbb1} that $A$
is a $J$-self-adjoint extension of $S$. Thus,
$A\in\Sigma_{J}^{\textsf{st}}$.

The last statement is an immediate consequence of
\cite[Proposition~1.4]{DM2}. This completes the proof.
\rule{2mm}{2mm}

Theorem~\ref{absd1} can be used for studying the spectral
properties of the operators $A\in\Sigma_{J}^{\textsf{st}}$. Recall
that if $S$ is simple and $\mu\in \rho(A_0)$, then it follows from
\eqref{bbb1} that for the components of the spectrum of $A=A_\cR$
one has
\begin{equation}\label{spec0}
 \mu\in\sigma_i(A) \quad\Leftrightarrow\quad
 0\in\sigma_i(M(\mu)-\cR)\quad (i=p,c,r),
\end{equation}
cf. \cite[Proposition~2.1]{DM}.

 \setcounter{equation}{0}
\section{The case of defect numbers $<2,2>$.}
Let ${J}$ and ${R}$ be a pair of \emph{anti-commuting} (i.e.,
$\{{J}, {R}\}=0$) fundamental symmetries in a Hilbert space
$\mathfrak{H}$ fundamental symmetries in a Hilbert space
$\mathfrak{H}$. Moreover, let $S$ be a densely defined symmetric
operator in $\mathfrak{H}$ with defect numbers $<2,2>$. In what
follows it is assumed that $S$ satisfies the following commutation
relations:
\begin{equation}\label{a2}
[S, {J}]=[S, {R}]=0.
\end{equation}

\subsection{ The descriptions of $\{C_\alpha\}_{\alpha\in{D}}$ and $\Sigma_{J}^{\textsf{st}}$.}
The first result describes the set $\{C_\alpha\}_{\alpha\in{D}}$ in
the present setting.

\begin{theorem}\label{uman4}
Let $S$ be a simple symmetric operator $S$ with defect numbers
$<2,2>$ which satisfies (\ref{a2}) and assume that the set
$\Upsilon$
\  is non-empty. Then the collection
$\{C_\alpha\}_{\alpha\in{D}}$  of all $C$-symmetries of $S$ in the
Krein space $(\mathfrak{H}, [\cdot,\cdot]_{{J}})$ coincides with the
set of operators ${C}_{\chi,\omega}$ defined by (\ref{sas44}).
\end{theorem}
\emph{Proof.} Using (\ref{neww12}) one can rewrite (\ref{sas44}) in
the form
\begin{equation}\label{uman3}
{C}_{\chi,\omega}=(\cosh\chi){J}+(\sinh\chi)(\cos\omega){JR}-i(\sinh\chi)(\sin\omega){R},
\end{equation}
where $\chi\in\mathbb{R}$ and $\omega\in[0,2\pi)$. It is known (see
\cite[Lemma 3.3]{AKD}) and easy to check that
${C}_{\chi,\omega}^2=I$ and ${J}{C}_{\chi,\omega}>0$. Furthermore,
it follows from (\ref{a2}) and \eqref{uman3} that $[S,
{C}_{\chi,\omega}]=0$. Therefore, $S$ has the
${C}_{\chi,\omega}$-symmetry for every choice of $\chi$ and
$\omega$.

Now, assume that $\Upsilon$ is non-empty. Then, by Lemma
\ref{neww41}, there exists a boundary triplet $(\mathcal{H},
\Gamma_0, \Gamma_1)$ of $S^*$ with the properties (\ref{neww50}).
This means (since $S$ has the ${C}_{\chi,\omega}$-symmetries) that
the operators ${C}_{\chi,\omega}$ have images ${\mathcal
C}_{\chi,\omega}$ in $\mathcal{H}$ determined by the formula
${\mathcal C}_{\chi,\omega}\Gamma_j=\Gamma_j{C}_{\chi,\omega}$
($j=0,1$). Considering this formula for $\omega=\pi/2$ and taking
the relations ${\mathcal J}\Gamma_j=\Gamma_j{J}$ and (\ref{uman3})
into account, one concludes that ${\mathcal R}\Gamma_j=\Gamma_j{R}$,
where ${\mathcal R}$ is a bounded operator in $\mathcal{H}$.
Therefore,
\begin{equation}\label{uman6}
{\mathcal C}_{\chi,\omega}=(\cosh\chi)\mathcal{J}+(\sinh\chi)(\cos\omega)\mathcal{JR}-i(\sinh\chi)(\sin\omega)\mathcal{R}
\end{equation}
and ${\mathcal C}_{\chi,\omega}^2=I$, $\mathcal{J}{\mathcal
C}_{\chi,\omega}>0$ in $\mathcal{H}$ (cf. Lemma \ref{neww67}).

Applying Lemma \ref{neww67} to ${\mathcal R}$, instead of ${\mathcal
J}$, it is seen that ${\mathcal R}$ is a fundamental symmetry in
$\mathcal{H}$. Moreover, $\{{\mathcal J}, {\mathcal R}\}=0$ since
$\{{J}, {R}\}=0$. Thus, ${\mathcal J}$ and ${\mathcal R}$ are
anti-commuting fundamental symmetries in ${\mathcal H}$.

Since $S$ has defect numbers $<2,2>$ the dimension of
${\mathcal H}$ is $2$.
Fix an orthonormal basis of ${\mathcal H}$ in which the matrix
representations of ${\mathcal J}$, ${\mathcal R}$, and $i{\mathcal
RJ}$ coincide with the Pauli matrices
\begin{equation}\label{pauli}
\sigma_3=\left(\begin{array}{cc} 1  & 0 \\
0 & -1
\end{array}\right), \quad \sigma_1=\left(\begin{array}{cc} 0  & 1 \\
1 & 0
\end{array}\right), \quad \sigma_2=\left(\begin{array}{cc} 0  & -i \\
i & 0
\end{array}\right),
\end{equation}
respectively. Then, according to 
\ (\ref{uman6}),
\begin{eqnarray}\label{sas8}
{\mathcal{C}}_{\chi,\omega}
=(\cosh\chi)\sigma_3+i(\sinh\chi)(\cos\omega)\sigma_2-i(\sinh\chi)(\sin\omega)\sigma_1
& & \nonumber\\ =\left(\begin{array}{cc}
\cosh\chi & (\sinh\chi)e^{-i\omega} \\
-(\sinh\chi)e^{i\omega} & -\cosh\chi
\end{array}\right). & &
\end{eqnarray}

Let $S$ have the $C$-symmetry in $(\mathfrak{H},
[\cdot,\cdot]_{{J}})$, i.e., $C\in\{C_\alpha\}_{\alpha\in{D}}$. Then
${\mathcal C}\Gamma_j=\Gamma_j{C}$ ($j=0,1$), where ${\mathcal
C}^2=I$ and ${\mathcal J}{\mathcal C}>0$. Considering ${\mathcal C}$
as $2\times{2}$-matrix we get
\begin{equation}\label{Cset1}
 {\mathcal C}^2=I \quad\text{and}\quad \sigma_3{\mathcal C}>0.
\end{equation}
It is easy to verify\footnote{see, e.g. the proof of
Lemma 3.5 in \cite{AKD}} that these conditions imposed on ${\mathcal
C}$ are equivalent to the presentation of ${\mathcal C}$ in the form
of (\ref{sas8}), i.e., $\mathcal{C}={\mathcal{C}}_{\chi,\omega}$ for
a certain choice of $\chi\in\mathbb{R}$ and $\omega\in[0,2\pi)$.
This means that ${C}\equiv{{C}}_{\chi,\omega}$, since the
correspondence ${{C}_\alpha}\rightarrow{\mathcal C}_\alpha$ is
injective by the simplicity of $S$; see Lemma \ref{neww41}.
Therefore, $\{C_\alpha\}_{\alpha\in{D}}=\{{C}_{\chi,\omega}\}$, where
$D=\mathbb{R}\times[0,2\pi)$ and $\alpha=(\chi,\omega)$.
Theorem \ref{uman4} is proved. \rule{2mm}{2mm}


\begin{corollary}\label{uman7}
Let $S$ satisfy the conditions of Theorem \ref{uman4}. Then
$(\mathcal{H}, \Gamma_0, \Gamma_1)$ is a boundary triplet for $S^*$
with properties (\ref{neww50}) if and only if the formulas
\begin{equation}\label{ea6a}
{\mathcal J}\Gamma_j:=\Gamma_j{J}, \qquad
{\mathcal R}\Gamma_j:=\Gamma_j{R}, \qquad j=0,1,
\end{equation}
correctly define intertwining operators ${\mathcal J}$ and ${\mathcal R}$ in $\mathcal{H}$.
\end{corollary}

\emph{Proof} This follows from the formula (\ref{uman3}) and
Theorem~\ref{uman4}. \rule{2mm}{2mm}

Recall that, if $S$ satisfies the conditions of Theorem \ref{uman4},
then there exist boundary triplets for $S^*$ with properties
(\ref{neww50}) or, what is equivalent, with the properties
(\ref{ea6a}) (Corollary \ref{uman7}). Fix one of them; $({\mathcal
H}, \Gamma_0, \Gamma_1)$. Then every $A\in\Sigma_{J}^{\textsf{st}}$
is determined by (\ref{as1}), where the corresponding ${\mathcal
J}$-{unitary operator} ${\mathcal K}$ is stable in the Krein space
$(\mathcal{H}, [\cdot,\cdot]_{{\mathcal J}})$ (Theorem
\ref{neww80}). Using Theorem~\ref{uman4} one can supplement the
results in Section 3 by giving explicit descriptions for the sets
 $\Sigma_{J}^{\textsf{st}}$ and $\Upsilon$.

In what follows the relationship $A\leftrightarrow{\mathcal K}$ is
indicated by using the notation $A_{\mathcal{K}}$ for
${J}$-self-adjoint operators $A$ determined by (\ref{as1}). The
operator ${\mathcal K}$ can be presented as a $2\times{2}$-matrix
$\mathcal{K}=(k_{ij})$. Since $\mathcal{K}$ is a ${\mathcal
J}$-unitary operator, the relation ${\mathcal J}={\mathcal
K}^*{\mathcal J}\mathcal{K}$ holds. Rewriting this in the matrix
form using the correspondence ${\mathcal J}\leftrightarrow\sigma_3$
(see (\ref{pauli}))
a simple analysis yields the following explicit formula for
$\mathcal{K}$:
\begin{equation}\label{as3a}
\mathcal{K}=\mathcal{K}(\zeta,
\phi, \xi, \omega)=e^{-i\xi}\left(\begin{array}{cc} -e^{-i\phi}\cosh\zeta & e^{-i\omega}\sinh\zeta \\
-e^{i\omega}\sinh\zeta & e^{i\phi}\cosh\zeta \end{array}\right),
\end{equation}
where $\zeta\in\mathbb{R}$, $\phi\in[0,\pi]$ and $\xi,
\omega\in[0,2\pi)$.

\begin{theorem}\label{p14}
Let the conditions of Theorem \ref{uman4} be satisfied. Then the
following statements are true:
\begin{itemize}

\item[(i)] if $A_\mathcal{K}\in\Sigma_{J}$, then its adjoint
$A_\mathcal{K}^*\in\Sigma_{J}$ is given by $ \mathcal{K}'(\zeta,
\phi, \xi, \omega)=\mathcal{K}(-\zeta, \phi, \xi, \omega)$;

\item[(ii)] $A_\mathcal{K}\in\Sigma_{J}$ is self-adjoint if and only if
$\zeta=0$, i.e.,
\begin{equation}\label{as3b}
\mathcal{K}=\mathcal{K}(0,
\phi, \xi, \omega)=e^{-i\xi}\left(\begin{array}{cc} -e^{-i\phi} & 0 \\
0 & e^{i\phi} \end{array}\right).
\end{equation}

\item[(iii)] $A_\mathcal{K}\in\Sigma_{J}$ belongs to $\Sigma_{J}^{\textsf{st}}\setminus\Upsilon$ if and
only if $\mathcal{K}={\mathcal{K}(\zeta, \phi, \xi, \omega)}$, where
\begin{equation}\label{sas55}
|\tanh\zeta|<|\cos\phi|;
\end{equation}

\item[(iv)] $A_\mathcal{K}\in\Sigma_{J}$ belongs to $\Upsilon$ if and only
if $\zeta=0$ and $\phi=\pi/2$, i.e.,
\begin{equation}\label{as3d}
\mathcal{K}=\mathcal{K}(0, {\pi}/{2}, \xi,
\omega)=\left(\begin{array}{cc}
{i}{e^{-i\xi}} & 0 \\
0 & {i}{e^{-i\xi}}
\end{array}\right).
\end{equation}
\end{itemize}
Furthermore, if the condition \eqref{sas55} in (iii) is satisfied,
then the operator $A_{\mathcal{K}(\zeta, \phi, \xi, \omega)}$ has
the ${C}_{\chi,\omega}$-symmetry, where the parameter $\chi$ is
(uniquely) determined by the equation
\begin{equation}\label{sas66}
{\cos\phi}\tanh\chi=-{\tanh\zeta}.
\end{equation}
\end{theorem}

\emph{Proof.} 
(i) This follows from (\ref{as3a}) by means of the identities
$A_\mathcal{K}^*=A_{{(\mathcal{K}^*)^{-1}}}$ and $(\cK^{-1})^*=JKJ$,
where $J=\sigma_3$; cf. \eqref{pauli}.

(ii) This immediately follows from (i).

(iii) It follows from Theorem~\ref{uman4} that $A\in\Sigma_{J}^{\textsf{st}}$ if and
only if $[A, {C}_{\chi,\omega}]=0$ for at least one
${C}_{\chi,\omega}$ defined by (\ref{uman3}).
Since $A=A_{\mathcal{K}}$ is defined by (\ref{as1}) (see Theorem \ref{neww80}), the corresponding operator
${\mathcal K}$ satisfies the commutation relation
$[{\mathcal{K}(\zeta, \phi, \xi, \omega)},
\mathcal{C}_{\chi,\omega}]=0$, where ${\mathcal C}_{\chi,\omega}\Gamma_j=\Gamma_j{C}_{\chi,\omega}$ ($j=0,1$).
A routine analysis of the last equality using (\ref{sas8}) and (\ref{as3a})
leads to the conclusion that $A_{\mathcal{K}}$ has the $C$-symmetry
($C\in\{{C_{\chi,\omega}}\}$) if and only if either $\zeta=0$,
$\phi=\frac{\pi}{2}$ (this case corresponds to the set $\Upsilon$;
cf. (iv)), or $|\tanh\zeta|<|\cos\phi|$. In the latter case the
operator $C$ can be chosen as ${C}_{\chi,\omega}$, where the
parameter $\chi$ is uniquely determined by the equation
(\ref{sas66}): this proves the last statement of the theorem.

(iv) By Theorem \ref{uman4} and the definition of $\Upsilon$ one has
$$
A\in\Upsilon \iff [A, {C}_{\chi,\omega}]=0, \quad \forall{\chi}\
\in\mathbb{R}, \ \omega\in[0,2\pi).
$$
The last relation is equivalent to $[\mathcal{K}, {\mathcal
J}]=[\mathcal{K}, {\mathcal R}]=0$ due to (\ref{uman6}). The first
condition $[\mathcal{K}, {\mathcal J}]=0$ imposed on the
${\mathcal{J}}$-unitary operator $\mathcal{K}$ means that
$\mathcal{K}$ is unitary. Hence, $A_\mathcal{K}$ is self-adjoint and
thus $\zeta=0$ and $\mathcal{K}$ is defined by (\ref{as3b}) (see
(ii)). Now, the second condition $[\mathcal{K}, {\mathcal R}]=0$
can be rewritten as $[\mathcal{K}(0, \phi, \xi, \omega),
\sigma_1]=0$ (since ${\mathcal R}\leftrightarrow\sigma_1$ by
(\ref{pauli})). Hence $-e^{-i\xi}e^{-i\phi}=e^{-i\xi}e^{i\phi}$ and
this equality holds if and only if $\phi=\frac{\pi}{2}$
($\phi\in[0,\pi]$).

The last statement was established while proving (iii). Hence,
Theorem \ref{p14} is proved. \rule{2mm}{2mm}

\subsection{Spectral analysis of $A\in\Sigma_{J}^{\textsf{st}}$.}
Let $(\mathcal{H}, \Gamma_0, \Gamma_1)$ be a boundary triplet for
$S^*$ with the properties (\ref{ea6a}) and let $M(\cdot)$ be the
corresponding Weyl function. Combining the definition (\ref{neww65})
of $M(\cdot)$ with (\ref{ea6a}) one arrives at the conclusion that
$[\mathcal{J}, M(\cdot)]=[\mathcal{R}, M(\cdot)]=0$. Passing to the
matrix representation $\mathcal{M}(\cdot)=(m_{ij}(\cdot))$ of
$M(\cdot)$ and using the correspondence ${\mathcal
J}\leftrightarrow\sigma_3$, ${\mathcal R}\leftrightarrow\sigma_1$
(see (\ref{pauli})), we get $[\sigma_3,
\mathcal{M}(\cdot)]=[\sigma_1, \mathcal{M}(\cdot)]=0$. This leads to
\begin{equation}\label{sas94}
M(\cdot)=m(\cdot)I,
\end{equation}
where $m(\cdot)$ is a scalar function defined on $\rho(A_0)$
($A_0=S^*\upharpoonright{\ker\Gamma_0}$).

Now consider an arbitrary ${C}\in\{{C}_{\chi,\omega}\}$ and the
corresponding decomposition
\begin{equation}\label{dede4}
{\mathfrak{H}}={\mathfrak
L}_+^{\chi,\omega}[\dot{+}]{\mathfrak
L}_-^{\chi,\omega}, \qquad {\mathfrak
L}_\pm^{\chi,\omega}=\frac{1}{2}(I\pm{C}_{\chi,\omega})\mathfrak{H}.
\end{equation}

Since $S$ and $S^*$ commute with ${C}_{\chi,\omega}$ they are
decomposed with respect to (\ref{dede4}):
\begin{equation}\label{sas77}
S=S_{+}(\chi,\omega)\dot{+}S_-(\chi,\omega), \qquad
S^*=S_{+}^{*}(\chi,\omega)\dot{+}S_-^{*}(\chi,\omega),
\end{equation}
where $S_{\pm}(\chi,\omega)=S\upharpoonright{{\mathfrak
L}_\pm^{\chi,\omega}}$ and $S_\pm^{*}(\chi,\omega)$ are adjoint of
the symmetric operators $S_{\pm}(\chi,\omega)$ acting in the
spaces\footnote{the spaces ${\mathfrak L}_{\pm}^{\chi,\omega}$ are
considered here with the original inner product $(\cdot,\cdot)$ on
$\mathfrak{H}$} ${\mathfrak L}_{\pm}^{\chi,\omega}$.

Let $A\in\Sigma_{J}^{\textsf{st}}$. Then $A=A_{\mathcal K}$ is
determined by (\ref{as3a}) and $A_{\mathcal
K}$ has the ${C}_{\chi,\omega}$-symmetry for a certain choice of
$\chi$ and $\omega$ (Theorem \ref{uman4}). Therefore,  $A_{\mathcal
K}$ is decomposed w.r.t. (\ref{dede4}):
\begin{equation}\label{uman10}
A_{\mathcal
K}=A_{\mathcal K}^+\dot{+}A_{\mathcal K}^-,
\end{equation}
where $S_{\pm}(\chi,\omega)\subset{A_{\mathcal
K}^\pm}\subset{S}_\pm^*(\chi,\omega)$. In this case, the ${\mathcal
J}$-unitary operator ${\mathcal{K}}$ is decomposed,
${\mathcal{K}}={\mathcal{K}_+}\dot{+}{\mathcal{K}_-}$, with respect
to the decomposition
\begin{equation}\label{neww10}
\mathcal{H}=\mathcal{H}_+^{\chi,\omega}[+]\mathcal{H}_-^{\chi,\omega}, \qquad
\mathcal{H}_\pm^{\chi,\omega}=\frac{1}{2}(I\pm{\mathcal{C}_{\chi,\omega}})\mathcal{H}
\end{equation}
of the Krein space $(\mathcal{H}, [\cdot,\cdot]_{{\mathcal J}})$
(cf. (\ref{dede4})). Since $\dim\mathcal{H}=2$, the subspaces
$\mathcal{H}_\pm^{\chi,\omega}$ are one-dimensional. Therefore,
${\mathcal{K}_\pm}=k_\pm{I}$ and eigenvalues $k_\pm$ of
$\mathcal{K}$ should satisfy the relations $
(\mathcal{K}-k_{\pm}I)(I\pm{\mathcal C}_{\chi,\omega})=0. $

A direct solution of the characteristic
equation $\det(\mathcal{K}-kI)=0$ gives
\begin{equation}\label{uman15}
k_{\pm}=e^{-i\xi}\left[\pm\sqrt{1-\sin^2\phi\cosh^2\zeta}+i\sin\phi\cosh\zeta\right].
\end{equation}

If, in particular, \eqref{sas55} is satisfied, a simple calculation
using (\ref{sas66}) leads to
\begin{equation}\label{as14}
k_+=-e^{-i\xi}e^{-it}, \, k_-=e^{-i\xi}e^{it}, \quad
e^{it}:=\frac{\cos\phi+i\sin\phi\cosh\chi}{|\cos\phi+i\sin\phi\cosh\chi|}.
\end{equation}
In this case the value of $\chi$ in (\ref{as14}) is uniquely
determined by $\zeta$ and $\phi$ (see (\ref{sas66})) and $t$ hence
can be considered as a function of $\zeta$ and $\phi$, i.e.,
$t=t(\zeta, \phi)(\in[0,2\pi))$.

On the other hand, if $\zeta=0$ and $\phi=\frac{\pi}{2}$ then in
\eqref{uman15} $k_\pm=ie^{-i\xi}$; cf. \eqref{as3d}. In this case
\eqref{as14} holds with $t=\frac{\pi}{2}$.

\begin{theorem}\label{wawa9}
Let the conditions of Theorem~\ref{uman4} be satisfied and assume
that $A_{\mathcal K}\in\Sigma_{{J}}^{\textsf{st}}$. Then the
spectrum of $A_{\mathcal K}$ is real and, moreover, $r\in\rho(A_0)$
belongs to the discrete spectrum of $A_{\mathcal K(\zeta, \phi, \xi,
\omega)}$ if and only if
\begin{equation}\label{uman9}
\left[\tan\frac{\xi+t}{2}+m(r)\right]\cdot\left[\cot\frac{\xi-t}{2}-m(r)\right]=0,
\end{equation}
where $m(\cdot)$ is given by (\ref{sas94}).

If, in particular, $A_{\mathcal K}\in
\Sigma_{{J}}^{\textsf{st}}\setminus\Upsilon$, then $t=t(\zeta,
\phi)$ is determined by (\ref{sas66}) and (\ref{as14}). Furthermore,
if $A_{\mathcal{K}}=A_{\mathcal{K}(0, {\pi}/{2}, \xi,
\omega)}\in\Upsilon$ and $\xi\not=\frac{\pi}{2}$, then
$r\in\rho(A_0)$ belongs to the discrete spectrum of
$A_{\mathcal{K}}$ if and only if
\begin{equation}\label{ad1}
\tan\frac{\xi+\frac{\pi}{2}}{2}+m(r)=0.
\end{equation}
In this case $A_{\mathcal{K}}$ coincides with self-adjoint operator
${A}_r$ defined by (\ref{sas99}).
\end{theorem}

\emph{Proof.} The reality of $\sigma(A_{\mathcal K})$ is a general
property of all ${J}$-self-adjoint operators with a ${C}$-symmetry
(see, e.g., Proposition \ref{sese1}).

Let $A_{\mathcal K}\in\Sigma_{{J}}^{\textsf{st}}$. Then $A_{\mathcal
K}$ admits the decomposition (\ref{uman10}) for certain $\chi$ and
$\omega$. In view of (\ref{as1}) and (\ref{sas77}), the
corresponding operators $A_{\mathcal K}^\pm$ are the restrictions of
${S}_\pm^*(\chi,\omega)$ onto
\begin{equation}\label{as1b}
\mathcal{D}(A_{\mathcal
K}^\pm)=\{f\in\mathcal{D}(S_\pm^*(\chi,\omega)) \ | \
k_\pm(\Gamma_1+i\Gamma_0)f=(\Gamma_1-i\Gamma_0)f\},
\end{equation}
cf. \eqref{as1bbb}. Rewriting the right-hand side of (\ref{as1b}) as
$i\frac{1+k_\pm}{1-k_\pm}\Gamma_0{f_\pm}=\Gamma_1{f_\pm}$ and taking
into account that
$$
i\frac{1+k_+}{1-k_+}=i\frac{1-e^{-i\xi}e^{-it}}{1+e^{-i\xi}e^{-it}}=-\tan\frac{\xi+t}{2},
\quad i\frac{1+k_-}{1-k_-}=\cot\frac{\xi-t}{2},
$$
where $t$ is given by (\ref{as14}), one gets
\begin{equation}\label{as1c}
\begin{array}{l}
\mathcal{D}(A_{\mathcal K}^+)=\{f\in\mathcal{D}(S_+^*(\chi,\omega))
\ | \ -\tan\frac{\xi+t}{2}\Gamma_0f=\Gamma_1f\}
\vspace{4mm}\\
\mathcal{D}(A_{\mathcal K}^-)=\{f\in\mathcal{D}(S_-^*(\chi,\omega))
\ | \ \cot\frac{\xi-t}{2}\Gamma_0f=\Gamma_1f\}.
\end{array}
\end{equation}

Using (\ref{dede4}), (\ref{neww10}) and recalling that
$\mathcal{C}_{\chi,\omega}\Gamma_j=\Gamma_j{C}_{\chi,\omega}$, it is
easy to see that the restrictions of the original boundary triplet
$(\mathcal{H}, \Gamma_0, \Gamma_1)$ onto the domains
$\mathcal{D}(S_\pm^*(\chi,\omega))$ give rise to the boundary
triplets $(\mathcal{H}_\pm^{\chi,\omega}, \Gamma_0, \Gamma_1)$ of
$S_\pm^*(\chi,\omega)$ in ${\mathfrak L}_\pm^{\chi,\omega}$.
Moreover, due to (\ref{sas94}), $m(\cdot)$ is the Weyl function of
$S_\pm(\chi,\omega)$ associated with the boundary triplets
$(\mathcal{H}_\pm^{\chi,\omega}, \Gamma_0, \Gamma_1)$. But, then by
Theorem~\ref{absd1} the formulas (\ref{as1c}) imply that
$r\in\rho(A_0)$ is an eigenvalue of $A_{\mathcal K}^+$ ($A_{\mathcal
K}^-$) if and only if $\tan\frac{\xi+t}{2}+m(r)=0$
($\cot\frac{\xi-t}{2}-m(r)=0$); see \eqref{spec0}. Now,
(\ref{uman9}) follows from the decomposition (\ref{uman10}).

The statement for $A_{\mathcal K}\in
\Sigma_{{J}}^{\textsf{st}}\setminus\Upsilon$ is clear. Assume that
$A_{\mathcal{K}}\in\Upsilon$. Then according to (\ref{as3d}) and
(\ref{uman15}) the eigenvalues $k_{\pm}$ of the operator
$\mathcal{K}$ coincide and $k_{\pm}=i{e^{-i\xi}}$. In particular,
$k_{\pm}\not=1$ precisely when $\xi\not=\frac{\pi}{2}$. Since
$\phi=\pi/2$ and $t=\pi/2$ in (\ref{as14}) one has
$-\tan\frac{\xi+\frac{\pi}{2}}{2}=\cot\frac{\xi-\frac{\pi}{2}}{2}$.
Therefore, (\ref{uman9}) reduces to (\ref{ad1}). In this case, the
(algebraic) multiplicity of the eigenvalue $r$ is $2$ and hence,
$A_{\mathcal{K}}=A_r$. Theorem~\ref{wawa9} is proved.
\rule{2mm}{2mm}

\begin{remark}\label{remm}
If $\xi=\frac{\pi}{2}$, then $\mathcal{K}(0, \frac{\pi}{2}, \frac{\pi}{2},
\omega)=I$ and the corresponding operator $A_{\mathcal{K}}$
coincides with $A_0=S^*\upharpoonright{\ker\Gamma_0}$; see
(\ref{as1bbb}). In this case, the formula (\ref{ad1}) means that
$m(r)$ has a pole at $r$, and this with a simply symmetric operator
$S$ means that $r\in\sigma(A_0)$.
\end{remark}

\subsection{Examples.} We start by describing a general procedure
which allows us to construct various examples illustrating the
results above. Our basic ingredients are: a symmetric operator $S_+$
with defect numbers $<1,1>$ acting in a Hilbert space
$\mathfrak{H}_+$; a boundary triplet $(\mathbb{C}, \Gamma_0^+,
\Gamma_1^+)$ of $S^*_+$; the Weyl function $m(\cdot)$ of $S_+$
associated with $(\mathbb{C}, \Gamma_0^+, \Gamma_1^+)$.

Let $\mathfrak{H}_-$ be a Hilbert space and let $W$ be a unitary map
of $\mathfrak{H}_-$ onto $\mathfrak{H}_+$. In the space
$\mathfrak{H}=\mathfrak{H}_+\oplus\mathfrak{H}$ consider the
operators
\begin{equation}\label{wawel1}
J=\left(\begin{array}{cc} I  &  0 \\
0 & -I
\end{array}\right), \quad R=\left(\begin{array}{cc} 0  &  W \\
W^{-1} &  0
\end{array}\right), \quad S=\left(\begin{array}{cc} S_+  &  0 \\
0 & {W^{-1}}S_+{W}
\end{array}\right).
\end{equation}

It is easy to see that $J$ and $R$ are anticommuting fundamental symmetries in $\mathfrak{H}$
and the symmetric operator $S$ satisfies the commutation relation (\ref{a2})
and it has defect numbers $<2,2>$ in $\mathfrak{H}$.

Let $S_+$ have real points of regular type. Then the operator $S$
also has real points of regular type. Hence, the set $\Upsilon$ is
non-empty (Theorem \ref{neww33}). Now one can use Corollary
\ref{uman7} to construct a boundary triplet of $S^*$ with the
properties (\ref{neww50}). More precisely, define the boundary
mappings $\Gamma_j^-:=\left(\begin{array}{c}
0 \\
\Gamma_j^+{R}
\end{array}\right)$. Then $(\mathbb{C}, \Gamma_0^-,
\Gamma_1^-)$ is a boundary triplet for $S^*_-=W^{-1}S^*_+{W}$ and
mappings $\Gamma_j=\left(\begin{array}{c}
\Gamma_j^+ \\
\Gamma_j^+{R}
\end{array}\right)
$ define a boundary triplet $({\mathbb{C}}^2, \Gamma_0, \Gamma_1)$
for $S^*$ which satisfies (\ref{ea6a}), where (cf. (\ref{pauli}))
\begin{equation}\label{ewa}
\mathcal{J}=\left(\begin{array}{cc}
1 & 0  \\
0 & -1
\end{array}\right)=\sigma_3, \qquad {\mathcal R}=\left(\begin{array}{cc}
0 &   1\\
1 & 0
\end{array}\right)=\sigma_1.
\end{equation}
By Corollary \ref{uman7}, the boundary
triplet $({\mathbb{C}}^2, \Gamma_0, \Gamma_1)$ satisfies (\ref{neww50}). Furthermore, the
Weyl function $M(\cdot)$ of $S$ associated with $(\mathcal{H}, \Gamma_0, \Gamma_1)$ is determined
by (\ref{sas94}), where $m(\cdot)$ is the Weyl function of $S_+$ associated with $(\mathcal{H}_+, \Gamma_0^+, \Gamma_1^+)$.

With these preparations, the spectral analysis of $J$-self-adjoint
operators with stable $C$-symmetries can be carried out by a
somewhat routine application of Theorem \ref{wawa9}. Observe, that
the corresponding spectral properties depend on the choice of the
initial symmetric operator $S_+$ (or, what is equivalent, on the
choice of the Weyl function $m(\cdot)$ of $S_+$).

The above considerations are illustrated with the following example.
\begin{example}\label{ex4}
Let $\mathfrak{H}_+=L_2^{\textsf{even}}(\mathbb{R})$ be the subspace
of even functions of $L_2(\mathbb{R})$ and define
$$
S_+=-\frac{d^2}{dx^2}, \qquad \mathcal{D}(S_+)=[{\stackrel{0}{W^2_2}}(\mathbb{R}_-)\oplus{\stackrel{0}{W^2_2}}(\mathbb{R}_+)]\cap{L_2^{\textsf{even}}(\mathbb{R})}.
$$

The adjoint $S_+^*=-\frac{d^2}{dx^2}$ has the domain
$\mathcal{D}(S_+^*)=W_2^{2}(\mathbb{R}\setminus\{0\})\cap{L_2^{\textsf{even}}(\mathbb{R})}$ and
$(\mathbb{C}, \Gamma_0^+, \Gamma_1^+)$ with
$$
\Gamma_0^+u(\cdot)=u(0),  \qquad \Gamma_1^+u(\cdot)=u'(+0)-u'(-0)=2u'(+0)
$$
defines a boundary triplet for $S_+^*$. With
$\mu\in\mathbb{C}\setminus\mathbb{R}$ the defect subspace
$\ker(S_+^*-\mu{I})$ coincides with the linear span of
$$
f_{\mu}=\left\{\begin{array}{cc}
 e^{i\tau{x}}, & x>0  \\
 e^{-i\tau{x}}, & x<0
 \end{array}\right. ,
$$ where $\tau=\sqrt{\mu}$ and $\textsf{Im}\ \tau>0.$ Since
$m(\mu)\Gamma_0^+f_{\mu}=\Gamma_1^+f_{\mu}$, the Weyl function of
$S_+$ associated with $(\mathbb{C}, \Gamma_0^+, \Gamma_1^+)$ is
given by
\begin{equation}\label{wawa96}
m(\mu)=2i\sqrt{\mu}.
\end{equation}

Let $\mathfrak{H}_-=L_2^{\textsf{odd}}(\mathbb{R})$ be the subspace
of odd functions of $L_2(\mathbb{R})$. According to (\ref{wawel1}),
the fundamental symmetry $J$ coincides with the space parity
operator $\mathcal{P}u(x)=u(-x)$ in
$\mathfrak{H}=L_2(\mathbb{R})=L_2^{\textsf{odd}}(\mathbb{R}){\oplus}L_2^{\textsf{even}}(\mathbb{R})$.
Choosing the unitary map $W :
L_2^{\textsf{odd}}(\mathbb{R}){\to}L_2^{\textsf{even}}(\mathbb{R})$
as $Wu=\textsf{sign}(x)u(x)$, one concludes that the fundamental
symmetry $R$ coincides with the multiplication by $\textsf{sign}(x)$
in $L_2(\mathbb{R})$. Now, the operator $S=-\frac{d^2}{dx^2}, \
\mathcal{D}(S)={\stackrel{0}{W^2_2}}(\mathbb{R}_-)\oplus{\stackrel{0}{W^2_2}}(\mathbb{R}_+)$
is symmetric in $L_2(\mathbb{R})$ and its adjoint
$S^*=-\frac{d^2}{dx^2}$ has the domain
$\mathcal{D}(S^*)=W_2^{2}(\mathbb{R}\setminus\{0\})$.

The boundary triplet $(\mathbb{C}^2, \Gamma_0, \Gamma_1)$ for $S^*$,
which is determined by
$$
\Gamma_0f=\Gamma_0(u+v)=\left(\begin{array}{c}
u(0) \\
v(+0)
\end{array}\right), \quad \Gamma_1f=2\left(\begin{array}{c}
u'(+0) \\
v'(0)
\end{array}\right),
$$
where $u$ and $v$ are, respectively, even and odd parts of $f$,
satisfies (\ref{ea6a}) with (\ref{ewa}). Moreover, all
$\mathcal{P}$-self-adjoint operators
$A_{\mathcal{K}}\in\Sigma_{\mathcal{P}}^{\textsf{st}}$ are
characterized by Theorem~\ref{p14} (statements (iii), (iv)).

The self-adjoint operator $A_0=S^*\upharpoonright{\ker\Gamma_0}$
coincides with the Friedrichs extension of $S$:
$$
A_0=-\frac{d^2}{dx^2}, \quad {\mathcal D}(A)=\{\,f(\cdot)\in{W_2^2}(\mathbb{R}\backslash\{0\}) :
 f(+0)=0, \ f(-0)=0\ \}.
$$

Applying Theorem \ref{wawa9} and taking the relations
$\sigma(A_0)=[0,\infty)$ and (\ref{wawa96}) into account, one
concludes that \emph{an arbitrary $A_{\mathcal K(\zeta, \phi, \xi,
\omega)}\in\Sigma_{\mathcal{P}}^{\textsf{st}}\setminus\Upsilon$ has
the essential spectrum on $[0,\infty)$ and a negative number $r$
belongs to the discrete spectrum of  $A_{\mathcal K}$ if and only
if}
$$
\left[\tan\frac{\xi+t}{2}-2\sqrt{|r|}\right]\cdot\left[\cot\frac{\xi-t}{2}+2\sqrt{|r|}\right]=0,
$$
where $t=t(\zeta, \phi)$ is determined by (\ref{sas66}) and
(\ref{as14}). The algebraic and the geometric multiplicities of $r$
are equal. 
\end{example}

\noindent \textbf{Acknowledgements.} The first author (S.H.) is
grateful for the support from the Finnish Cultural Foundation, South
Ostrobothnia Regional fund. The second author (S.K.) expresses his
gratitude to the Academy of Finland (projects 128059, 132533) and
project JRP IZ73Z0 (28135) of SCOPES 2009-2012 for the support.

\end{document}